\input amstex
\documentstyle{mtn7-am}
\pageno=1

\loadbold 


\define\0{\bold0}
\define\1{\bold1}
\define\sB{\script B}
\define\sP{\script P}
\define\Map{\operatorname{Map}}
\define\lip{\operatorname{lip}}
\define\blip{\operatorname{\hbox{\bf lip}}}
\define\Hom{\operatorname{Hom}}
\define\Low{\operatorname{Low}}
\define\Up{\operatorname{Up}}
\define\epi{\operatorname{epi}}
\define\Conv{\operatorname{Conv}}
\define\Conc{\operatorname{Conc}}

\define\Sect#1#2#3{\par
\smallskip
\parshape=2 1mm 126mm 1mm 126mm
\ninepoint\noindent{\bf#1.}\enspace#2\dotfill\dotfill
\noindent\vskip-\normalbaselineskip\parshape=1 0mm 134mm {}
\hfill\ninepoint#3}
\define\SubSect#1#2#3{\par
\smallskip
\parshape=2 5mm 122mm 5mm 122mm
\ninepoint\noindent{\bf#1.}\enspace#2\dotfill\dotfill
\noindent\vskip-\normalbaselineskip\parshape=1 0mm 134mm {}
\hfill\ninepoint#3}
\define\RefSect#1{\par
\smallskip
\parshape=1 1mm 126mm
\ninepoint\noindent Список литературы\dotfill\dotfill
\noindent\vskip-\normalbaselineskip\parshape=1 0mm 134mm {}
\hfill\ninepoint#1}
\let\emptyset\varnothing
\let\phi\varphi
\let\epsilon\varepsilon
\let\wh\widehat
\let\wt\widetilde
\let\ol\overline
\let\<\langle
\let\>\rangle

\def\Russpages{758--797}
\def\Russvol{69}
\def\Russno{5}
\def\Englno{56}

\def\Russyr{2001}
\def\Englyr{2001}
\def\copyfee{25.00}

\comment
\endcomment
\titleinfo
{Idempotent functional analysis:}
\extratitleline
{An algebraic approach}

\def\righttext{IDEMPOTENT FUNCTIONAL ANALYSIS}
\def\lefttext{G.~L.~LITVINOV \ {\sl et al.}}
\def\shortauthors{Litvinov, Maslov, Shpiz}

\authorinfo
{G.~L.~Litvinov, V.~P.~Maslov, and G.~B.~Shpiz}
{Received May~17, 2000%
}

\topmatter
\abstract\nofrills
\ninepoint
{\bf Abstract}---
This paper is devoted to Idempotent Functional Analysis, which is an
``abstract'' version of Idempotent Analysis developed by V.~P.~Maslov and
his collaborators. We give a brief survey of the basic ideas of Idempotent
Analysis.  The correspondence between concepts and theorems of the
traditional Functional Analysis and its idempotent version is discussed in
the spirit of N.~Bohr's correspondence principle in Quantum Theory.
We present an algebraic approach to Idempotent Functional Analysis.
Basic notions and results are formulated in algebraic terms; the essential
point is that the operation of idempotent addition can be defined for
arbitrary infinite sets of summands. We study idempotent analogs of the
basic principles of linear functional analysis and results on the
general form of a linear functional and scalar products in idempotent
spaces.

\medskip\noindent
Published in {\it Math.\ Notes}, vol.~69, no.~5 (2001), p.~696--729.

\comment
 Ключевые слова: {\it идемпотентный анализ, функциональный
анализ, линейный функционал, скалярное произведение.}
\endcomment
\par\removelastskip\medskip\vskip1ex\noindent
{\eightsmc Key words}: \it 
Idempotent Analysis, Functional Analysis, linear functional, scalar
product.
\endabstract
\endtopmatter

\document
\par\removelastskip\vskip18pt plus3pt minus3pt


\rightline{\ninepoint And above all I value Analogies,}
\vskip-1.5pt
\rightline{\ninepoint my most reliable teachers.}
\rightline{\ninepoint\it J.~Kepler}

\ahead
 Table of contents
\endahead
\bgroup
\Sect{1}{Introduction}{2}
\SubSect{1.1}{Heuristics}{2}
\SubSect{1.2}{Idempotent Analysis}{3}
\SubSect{1.3}{The superposition principle and linear problems}{5}
\SubSect{1.4}{Idempotent Functional Analysis}{7}
\Sect{2}{Idempotent semigroups and partial order}{8}
\Sect{3}{Idempotent semirings, semifields, and quasifields}{16}
\SubSect{3.1}{Idempotent semirings and semifields}{16}
\SubSect{3.2}{Examples}{17}
\SubSect{3.3}{Complete semirings}{17}
\SubSect{3.4}{Quasifields}{18}
\SubSect{3.5}{Completion of semirings}{19}
\SubSect{3.6}{Examples}{20}
\Sect{4}{Idempotent semimodules and spaces}{20}
\SubSect{4.1}{Basic definitions}{20}
\SubSect{4.2}{Complete and standard semimodules}{21}
\SubSect{4.3}{Idempotent spaces}{22}
\SubSect{4.4}{Linear maps and functionals}{22}
\SubSect{4.5}{Idempotent semimodules and spaces
associated with vector lattices}{23}
\SubSect{4.6}{Examples}{24}
\Sect{5}{The structure of $a$-linear functionals on idempotent spaces}{25}
\SubSect{5.1}{The basic construction}{25}
\SubSect{5.2}{The basic theorem on the structure of functionals}{26}
\SubSect{5.3}{Theorems of the Hahn--Banach type}{27}
\SubSect{5.4}{Analogs of the Banach--Steinhaus theorem
and the closed graph theorem}{28}
\SubSect{5.5}{The scalar product}{29}
\SubSect{5.6}{The skew-scalar product and elements of duality}{30}
\SubSect{5.7}{Examples}{31}
\Sect{6}{Commentary}{32}

\egroup

\ahead
1.
 INTRODUCTION
\endahead

\bhead
1.1.
 Heuristics
\endbhead
Idempotent Functional Analysis is an abstract version of Idempotent
Analysis in the sense of~\cite{1--8}. Idempotent Analysis is closely
related to optimal control theory, optimization theory, convex analysis,
the theory of vector lattices and ordered algebraic systems.
Ideas, results, and terminology accumulated in these fields are taken
into account where possible, but Idempotent Functional Analysis is considered
here as a part of Idempotent Mathematics (see below) and presented
in the framework of a different paradigm.

This paradigm is expressed by the {\it correspondence principle\/}
\cite{9, 10}. This principle is similar to the well-known correspondence
principle of N.~Bohr in Quantum Theory (and closely related to it).
Actually, there exists a heuristic correspondence between
important, interesting, and useful constructions and results of
traditional mathematics over fields and analogous constructions and results
over idempotent semirings and semifields (i.e., semirings and semifields
with idempotent addition; for rigorous definitions see below).

A systematic and consistent application of the ``idempotent''
correspondence principle leads to a variety of results, often quite
unexpected. As a result, in parallel with traditional mathematics over
rings, its ``shadow'', Idempotent Mathematics, appears. This ``shadow''
stands approximately in the same relation to traditional mathematics as
classical physics to Quantum Theory. In many respects Idempotent
Mathematics is simpler than traditional mathematics. However, the transition
from traditional concepts and results to their idempotent analogs is often
nontrivial. A correct formulation of these analogs is sometimes the most
difficult step. In this sense, Idempotent Mathematics resembles mathematics
over $p$-adic fields.

Let~$\Bbb R$ be the field of real numbers and $\Bbb R_+$ the semiring of all
nonnegative real numbers (with respect to the ordinary sum and product).
The change of variables $x\mapsto u=h\,\ln x$
defines a map
$\Phi_h\:\Bbb R_+\to S=\Bbb R\cup\{-\infty\}$.
Let the addition
and multiplication operations be mapped from $\Bbb R$ to~$S$ by~$\Phi_h$,
i.e., let
$$
u\oplus_hv=h\,\ln(\exp(u/h)+\exp(v/h)),\quad
u\odot v=u+ v,\quad
\0=-\infty=\Phi_h(0),\quad
\1=0=\Phi_h(1).
$$
It can
easily be checked that $u\oplus_hv\to\max\{u,v\}$ as $h\to 0$, and $S$
forms a semiring with respect to the addition $u\oplus v=\max\{u,v\}$ and
the multiplication $u\odot v=u+ v$ with zero $\0= -\infty$ and unit $\1=
0$.  Denote this semiring by~$\Bbb R_{\max}$; it is {\it idempotent}, i.e.,
$u\oplus u= u$ for all its elements.\footnote{The semiring $\Bbb R_{\max}$ is
actually a semifield; also, it is often called the $(\max,+)$ algebra.}
The analogy with quantization is obvious; the parameter~$h$ plays the
r\^{o}le of the Planck constant, so~$\Bbb R_+$ (or~$\Bbb R$)
can be viewed as
a ``quantum object'' and~$\Bbb R_{\max}$ as the result of its ``dequantization''.
A similar procedure gives the semiring $\Bbb R_{\min}=\Bbb R\cup\{+ \infty\}$
with the operations $\oplus= \min$, $\odot=+$;
in this case $\0=+ \infty$, $\1= 0$.
The semirings~$\Bbb R_{\max}$ and~$\Bbb R_{\min}$
are isomorphic.
Connections
with physics and imaginary values of the Planck constant are discussed
below in Sec.~1.3.
The idempotent semiring $\Bbb R\cup\{-\infty\}\cup\{+ \infty\}$
with the operations $\oplus= \max$, $\odot= \min$ is
obtained as a result of a ``second dequantization'' of $\Bbb R$ (or~$\Bbb R_+$).
  Dozens of interesting examples of nonisomorphic idempotent
semirings may be cited, as well as a number of standard methods of deriving
new semirings from them (see~\cite{1--9, 11--21} and below).

\bhead
1.2.
 Idempotent Analysis
\endbhead
Let~$S$ be an arbitrary semiring with idempotent addition~$\oplus$ (which
is always assumed to be commutative), multiplication~$\odot$, zero~$\0$, and unit~ $\1$.
The set~$S$ is supplied with the {\it standard partial
order\/}~$\preccurlyeq$: by definition, $a\preccurlyeq b$,
if (and only if) $a\oplus b= b$.
Thus all elements of~$S$ are positive: $\0\preccurlyeq a$ for all $a\in S$.
 Due
to the existence of this order, Idempotent Analysis is closely related to
lattice theory~\cite{22, 23}, the theory of vector lattices, and the theory
of ordered spaces~\cite{21--27}.
 Moreover, this partial order allows to
model a number of basic notions and results of Idempotent Analysis at
the purely algebraic level, see~\cite{11--14}; in this paper we develop
this line of reasoning systematically.

Calculus deals mainly with functions whose values are numbers.
 The
idempotent analog of a numerical function is a map $X\to S$, where~$X$ is an arbitrary set and~$S$ is an idempotent semiring.
 Functions with values
in~$S$ can be added, multiplied by each other, and multiplied by elements
of~$S$ pointwise.

The idempotent analog of a linear functional space is a set of~$S$-valued
functions that is closed under addition of functions and multiplication of
functions by elements of~$S$, or an~$S$-semimodule.
 Consider, e.g., the~$S$-semimodule $\sB(X,S)$ of functions $X\to S$ that are bounded in
the sense of the standard order on~$S$.

If $S= \Bbb R_{\max}$, then the idempotent analog of integration is defined by the
formula
$$
I(\phi)
=\int_X^\oplus\phi(x)\,dx
=\sup_{x\in X}\phi(x),
\tag1.1
$$
where $\phi\in\sB(X,S)$. Indeed, a Riemann sum of the form
$$
\sum_i\phi(x_i)\cdot\sigma_i
$$
corresponds to the expression
$$
\bigoplus_i\phi(x_i)\odot\sigma_i
=\max_i\{\phi(x_i)+\sigma_i\},
$$
which tends to the right-hand side of~\thetag{1.1} as $\sigma_i\to 0$. Of
course, this is a purely heuristic argument.

Formula~\thetag{1.1} defines the idempotent integral not only for functions
taking values in~$\Bbb R_{\max}$, but also in the general case when any
bounded (from above) subset of~$S$ has a least upper bound; semirings of
this type are called {\it boundedly complete}.

An idempotent measure on~$X$ is defined by
$$
m_\psi(Y)=\sup_{x\in Y}\psi(x),
$$
where $\psi\in B(X,S)$.
 The integral with respect to this
measure is defined by
$$
I_\psi(\phi)
=\int^\oplus_X\phi(x)\,dm_\psi
=\int_X^\oplus\phi(x)\odot\psi(x)\,dx
=\sup_{x\in X}\bigl(\phi(x)\odot\psi(x)\bigr).
\tag1.2
$$

Obviously, if $S= \Bbb R_{\min}$,
then the standard order~$\preccurlyeq$
is opposite to
the usual order~$\le$, so in this case Eq.~\thetag{1.2} assumes the form
$$
\int^\oplus_X\phi(x)\,dm_\psi
=\int_X^\oplus\phi(x)\odot\psi(x)\,dx
=\inf_{x\in X}\bigl(\phi(x)\odot\psi(x)\bigr),
\tag1.3
$$
where~$\inf$ is understood in the sense of the usual order~$\le$.

The functionals $I(\phi)$ and $I_\psi(\phi)$ are linear over~$S$;
their values correspond to limits of Lebesgue (or Riemann) sums.
The formula for $I_\psi(\phi)$ defines the
idempotent scalar product of the functions~$\psi$ and~$\phi$.
 Various
idempotent functional spaces and an idempotent version of the theory of
distributions can be constructed on the basis of idempotent
integration~\cite{1--8}.
 The analogy between idempotent and probability measures leads to
spectacular parallels between optimization theory and probability theory.
 For
example, the Chapman--Kolmogorov equation corresponds to the Bellman
equation (see the survey of Del~Moral~\cite{28} and~\cite{29--35, 12, 13}).
Many other idempotent analogs may be cited (in particular, for the basic
constructions and theorems of functional analysis).
 For instance, the
Legendre transform is nothing but the $\Bbb R_{\max}$ version of the Fourier
transform, see~\cite{1--8}.

Indeed, suppose $S= \Bbb R_{\max}$, $G= \Bbb R^n$;
let~$G$ have a topological group structure.
 The ordinary Fourier--Laplace transform is defined as
$$
\phi(x)\mapsto\wt\phi(\xi)
=\int_Ge^{i\xi\cdot x}\phi(x)\,dx,
\tag1.4
$$
where~$e^{i\xi\cdot x}$ is a character of the group~$G$, i.e., a solution
of the following functional equation:
$$
f(x+y)=f(x)f(y).
$$
The idempotent analog of this equation is
$$
f(x+y)=f(x)\odot f(y)=f(x)+f(y),
$$
so ``continuous idempotent characters'' are linear functionals of the
form
$$
x\mapsto\xi\cdot x={\xi_1x_1+\dots+\xi_nx_n}.
$$
 As a result,
the transform in~\thetag{1.4} assumes the form
$$
\phi(x)\mapsto\wt\phi(\xi)
=\int_G^\oplus\xi\cdot x\odot\phi(x)\,dx
=\sup_{x\in G}\bigl(\xi\cdot x+\phi(x)\bigr).
\tag1.5
$$
The transform in~\thetag{1.5} is nothing but the {\it Legendre transform\/}
(up to some notation); transforms of this kind establish a
correspondence between the Lagrangian and the Hamiltonian formulations
of classical mechanics.

\bhead
1.3.
The superposition principle and linear problems
\endbhead
 The basic equations of Quantum Theory
are linear (the superposition principle).
 The Hamilton--Jacobi equation, the basic equation of classical
mechanics, is nonlinear in the usual sense.
 However it is linear over the
semiring $\Bbb R_{\min}$.
 Also, different versions of the Bellman equation, the basic
equation of optimization theory, are linear over suitable idempotent
semirings (V.~P.~Maslov's idempotent superposition
principle), compare~\cite{1--8}.
 For instance, the finite-dimensional
stationary Bellman equation can be written in the the form $X=H\odot X\oplus F$,
where $X,H,F$ are matrices with coefficients in
an idempotent semiring~$S$ and the unknown matrix~$X$ is determined by~$H$
and~$F$~\cite{16}.
 In particular, standard problems of dynamic programming
and the well-known shortest path problem correspond to the cases $S= \Bbb R_{\max}$ and $S= \Bbb R_{\min}$, respectively.
 In~\cite{16}, it was shown that the main optimization
algorithms for finite graphs correspond to standard methods for solving
systems of linear equations of this type (i.e., over semirings).
Specifically, Bellman's shortest path algorithm
corresponds to a version of Jacobi's algorithm, Ford's algorithm
corresponds to the Gauss--Seidel iterative scheme, etc.; see also~\cite{1,
5, 7, 8, 16--21, 12, 13, 74}.

The linearity of the Hamilton--Jacobi equation over $\Bbb R_{\min}$
(and~$\Bbb R_{\max}$) is closely related to the (ordinary) linearity of
the Schr\"{o}dinger equation.  Consider the classical dynamical system
specified by the Hamiltonian
$$
H=H(p,x)=\sum_{i=1}^N\frac{p^2_i}{2m_i}+V(x),
\tag1.6
$$
where $x=(x_1,\dots,x_N)$ are generalized coordinates, $p=(p_1,\dots,p_N)$ are generalized momenta, $m_i$ are generalized masses, and
$V(x)$ is the potential.
 In this case the Lagrangian $L(x,\dot x,t)$ has the form
$$
L(x,\dot x,t)
=\sum^N_{i=1}m_i\frac{\dot x_i^2}2-V(x),
\tag1.7
$$
where $\dot x=(\dot x_1,\dots,\dot x_N)$, $\dot x=dx/dt$.
 The action functional $S(x,t)$ has the form
$$
S(x,t)=\int^t_{t_0}L(x,\dot x,t)\,dt,
\tag1.8
$$
where the integration is performed along a trajectory of the system.
  The
classical equations of motion are derived as the stationarity conditions
for the action functional (the Hamilton principle, or least action
principle), see, e.g.,~\cite{37}.

The action functional can be regarded as a function taking the set of
curves (trajectories) to the set of real numbers.
 Assume that its range
lies in the semiring~$\Bbb R_{\min}$.
 In this case the minimum of the action
functional can be viewed as the idempotent integral of this function over
the set of trajectories or the idempotent analog of the Feynman path
integral.
 Thus the least action principle can be considered as the idempotent
version of the well-known Feynman approach to quantum mechanics (which is
presented, e.g., in~\cite{38}); here, one should remember that the
exponential function involved in the Feynman integral is monotone on the
real axis.
 The representation of a solution to the Schr\"{o}dinger
equation in terms of the Feynman integral corresponds to the
Lax--Ole\u{\i}nik formula for solving the Hamilton--Jacobi equation
(see below).

Since
$$
\frac{\partial S}{\partial x_i}=p_i,
\qquad
\frac{\partial S}{\partial t}=-H(p,x),
$$
the following Hamilton--Jacobi equation holds:
$$
\frac{\partial S}{\partial t}
+H\biggl(\frac{\partial S}{\partial x_i},x_i\biggr)
=0.
\tag1.9
$$

Quantization (see, e.g.,~\cite{39}) leads to the Schr\"{o}dinger equation
$$
-\frac\hbar i\frac{\partial\psi}{\partial t}
=\wh H\psi=H(\wh p_i,\wh x_i)\psi,
\tag1.10
$$
where $\psi=\psi(x,t)$ is the wave function, i.e., a time-dependent
element of the Hilbert space $L^2(\Bbb R^N)$, and
$\wh H$~is the energy
operator obtained by substituting the momentum operators
$$
\wh p_i=\frac\hbar i\frac\partial{\partial x_i}
$$
and the coordinate
operators $\wh x_i\:\psi\mapsto x_i\psi$ for the variables
$p_i$ and $x_i$ in the Hamiltonian function~$H(p,x)$, respectively.
 This equation is
linear in the ordinary sense (the quantum superposition principle).
 The
standard procedure of limit transition from the Schr\"{o}dinger equation to
the Hamilton--Jacobi equation is to use the following ansatz for the wave
function:  $\psi=a(x,t)e^{(i/\hbar)S(x,t)}$, and to expand Eq.~\thetag{1.10} in powers of $\hbar \to 0$ (the `semiclassical' limit).

Instead of doing this, we switch to imaginary values of the Planck
constant~$\hbar$ by the substitution $h = i\hbar$, assuming $h > 0$.
 Thus the Schr\"{o}dinger equation~\thetag{1.10} turns to an analog of the
heat equation:
$$
h\frac{\partial u}{\partial t}
=H\biggl(-h\frac\partial{\partial x_i},\wh x_i\biggr)u,
\tag1.11
$$
where the real-valued function~$u$ corresponds to the wave function~$\psi$.
 A similar idea (the switch to imaginary time) is used in
Euclidean quantum field theory (see, e.g.,~\cite{40, 41}); let us remember
that time and energy are dual quantities.

The linearity of Eq.~\thetag{1.10} implies that of Eq.~\thetag{1.11}.
 Thus if~$u_1$ and~$u_2$ are solutions of~\thetag{1.11}, then so is their linear
combination
$$
u=\lambda_1u_1+\lambda_2u_2.
\tag1.12
$$

Let $S=-h\,\ln u$ or $u = e^{-S/h}$ as in Sec.~1.1 above.
 It can easily
be checked that Eq.~\thetag{1.11} thus assumes the form
$$
\frac{\partial S}{\partial t}
=V(x)+\sum^N_{i=1}\frac1{2m_i}
\biggl(\frac{\partial S}{\partial x_i}\biggr)^2
-h\sum^n_{i=1}\frac1{2m_i}
\frac{\partial^2S}{\partial x^2_i}.
\tag1.13
$$
This equation is nonlinear in the ordinary sense.
 However, if~$S_1$ and~$S_2$
are its solutions, then so is the function
$$
S=\lambda_1\odot S_1\oplus_h\lambda_2\odot S_2,
\tag1.14
$$
obtained from~\thetag{1.12} by means of our substitution
$S=-h\,\ln u$.
  Here the
generalized multiplication~$\odot$ coincides with ordinary addition and
the generalized addition is the image of ordinary addition
under the above change of variables.
  As $h \to 0$, we obtain the
operations of the idempotent semiring $\Bbb R_{\min}$, i.e.,
$\odot=+$, $\oplus= \min$, and
Eq.~\thetag{1.13}
turns to the Hamilton--Jacobi equation~\thetag{1.9},
since the third term in the right-hand side of
Eq.~\thetag{1.13}
vanishes.

Thus it is natural to consider the limit function $S=\lambda_1\odot S_1\oplus\lambda_1\odot S_2$
as a solution of the Hamilton--Jacobi equation and
to expect that this equation may be treated as linear over~$\Bbb R_{\min}$.
 This
argument (clearly, a heuristic one) can be extended to equations of a more
general form.
 For a rigorous treatment of (semiring) linearity for these
equations, see~\cite{6--8} and also~\cite{3}.
 Notice that if~$h$ is changed
to~$-h$, then the resulting Hamilton--Jacobi equation is linear over
$\Bbb R_{\max}$.

The idempotent superposition principle indicates that there exist important
problems that are linear over idempotent semirings.

\bhead
1.4.
 Idempotent functional analysis
\endbhead
Idempotent Functional Analysis is an analog of traditional Functional
Analysis in the framework of Idempotent Mathematics.
 Here we formulate a number
of well-known results of Idempotent Analysis at a new, more abstract level;
also, we consider problems that were not discussed in the earlier
literature on Idempotent Analysis.

The most important results of Idempotent Analysis obtained so far are fixed
point theorems and the spectral theory of linear operators on specific
idempotent semimodules.
 These operators satisfy additional conditions of
continuity type (see, e.g.,~\cite{2--8}) or of algebraic regularity
type in the spirit of~\cite{14, 6, 11--13} and the theory of~$C^*$
and~$W^*$ algebras.
  The known theorems on spectra of operators depend on the
idempotent analog of the integral representation of operators (an analog of
the L.~Schwartz kernel theorem), which in turn employs the integral
representation of linear functionals.

At the moment, the best studied idempotent semimodules are semimodules of
bounded functions (with values in idempotent semirings), idempotent
semimodules of continuous or semicontinuous functions (taking values in the
extended real axis and other idempotent semirings) that are defined on open
subsets of compacta with various boundary conditions, and free
finite-dimensional idempotent semimodules (the idempotent linear algebra).
In addition, a study of a class of semimodules that can be regarded as
idempotent analogs of Sobolev spaces was carried out in~\cite{42, 43, 6} in
connection with the study of the Hamilton--Jacobi equations.
 For these
semimodules, neither theorems on the general form of functionals and
integral representations of operators nor spectral theorems are proved, but
existence and uniqueness theorems for fixed points of linear operators were
obtained under some additional conditions.

Note that, as a rule, subsemimodules of semimodules of the classes studied
so far do not belong to the same classes (for instance, a subsemimodule of
a free semimodule generally is not free).
 Another important example that cannot be
treated by known methods is the semimodule of integrable
functions on a measure space.
 In this case any (idempotent) linear
functional that is continuous (with respect to the $L^1$~norm) is constant.

In this paper we present an algebraic approach to Idempotent Functional
Analysis: basic notions and results are ``simulated'' in algebraic terms.
Elements of this approach can be traced back to~\cite{11--13} and
especially~\cite{14}.
 The essential point is that the operation of
idempotent addition can be defined for an infinite set of summands if the
semiring is complete as an ordered set.
  In this case, the continuity
property of functionals and maps can be simulated by the preservation of
infinite sums under a suitable completion.
  Central to our
analysis is a natural class of abstract semimodules that are idempotent
analogs of vector spaces.
 We call semimodules of this kind {\it idempotent
spaces\/} (see Sec.~4 below).
 The class of idempotent spaces contains
most of the examples important for Functional Analysis, and subspaces of
idempotent spaces belong to the same class.
  In particular, idempotent
subspaces of (topological) vector lattices provide a very important example
of idempotent spaces.

Further, in this paper we present idempotent versions of the basic results
concerning linear functionals and scalar products, including the theorem on
the general form of a linear functional and idempotent analogs of the
Hahn--Banach and Riesz--Fischer theorems.
 We also present analogs of the
Banach--Steinhaus and the closed graph theorems.
 In forthcoming papers we
will apply the ``topological'' approach in the spirit of~\cite{1--8} and
construct abstract Idempotent Functional Analysis starting from basic
notions and results and leading up to analogs of A.~Grothendieck's results
on topological tensor products, kernel spaces, and operators.

For additional comments and historical remarks see Sec.~6 (Commentary).

\ahead
2.
 IDEMPOTENT SEMIGROUPS AND PARTIAL ORDER
\endahead

\chead
2.1
\endchead
Recall that a {\it semigroup\/} is a nonempty set endowed with an
associative operation called addition (in additive semigroups) or
multiplication (in multiplicative semigroups).
 A semigroup with a neutral
element (called its zero or unit, respectively) is called a {\it
monoid\/}; see, e.g.,~\cite{22}.

\definition{Definition~2.1}
An {\it idempotent semigroup\/} is an additive semigroup~$S$ with
commutative addition~$\oplus$ such that $x \oplus x = x$ for all $x \in S$.
If this semigroup is a monoid, then its neutral element is denoted by~$\0$
or~$\0_S$.
\enddefinition

For a rich collection of examples of idempotent semigroups see Sec.~2.9
below.


\chead
2.2
\endchead
Any idempotent semigroup is a partially ordered set with respect to the
standard order defined below.

\definition{Definition~2.2}
The {\it standard order\/} on an idempotent semigroup~$S$ is the partial
order~$\preccurlyeq$
such that $x\preccurlyeq y$
if and only if $x\oplus y= y$.
\enddefinition

From now on, we assume that all idempotent semigroups (and semirings) are
ordered in this way and the notation~$\prec$, $\succ$, $\succcurlyeq$ has
the obvious meaning.
 For instance, the relation $x\prec y$ means that
$x\preccurlyeq y$ and $x\ne y$.

For basic definitions and results of the theory of ordered sets see, e.g.,
\cite{23}.
 It can easily be checked that Definition~2.2 is self-consistent
and defines an order relation.
  It follows from this definition that $\0\preccurlyeq x$ for all $x\in S$ if~$S$ is a monoid.
 By~$\Bbb I$ or~$\Bbb I_S$
denote the (unique) element of~$S$
such that
$x\preccurlyeq\Bbb I$ for all $x\in S$
(if this element exists).
 It is clear that
$\0=\inf S=\sup\emptyset$ and
$\Bbb I=\sup S=\inf\emptyset$,
where~$\emptyset$ is the empty subset of~$S$ and
$\inf$~and~$\sup$ denote the greatest lower bound and the least upper bound, respectively.

Any idempotent semigroup~$S$ is a~$\vee$-semilattice (or upper
semilattice), i.e., for any $x, y \in S$ the set $\{x, y\}$ has the least
upper bound $\sup \{x, y\} = x \vee y$.
 Obviously,
$$
\sup\{x,y\}=x\oplus y.
$$
Thus the class of all idempotent semigroups coincides with the class of all
upper semilattices.

\remark{Remark 2.1}
Although formally the notions of idempotent semigroup and semilattice
coincide, classical lattice theory is aimed at the study of lattices of
subsets (subspaces, ideals, etc\.).
 Idempotent Analysis deals primarily
with sublattices of ordered vector and functional spaces.
 These two
approaches are in a relation similar to that of measure theory and
probability theory: the object is common but the sources of interesting
problems are different.
\endremark

\chead
2.3
\endchead
Let an idempotent semigroup~$S$ be a lattice with respect to the standard
order, i.e., let a greatest lower bound $\inf \{x, y\}$, denoted by $x
\wedge y$, be defined for any two elements $x, y \in S$ along with their
least upper bound $x \oplus y$.
 In this case we call the semigroup~$S$ a {\it lattice
semigroup\/} or simply (by abuse of terminology) a {\it
lattice\/}.

\definition{Definition~2.3}
{\it An idempotent semigroup\/}~$S^\circ$ is called {\it ordinally dual\/} (or
$o$-{\it dual\/}) {\it to a lattice semigroup\/}~$S$ if it consists of the set of all
elements of~$S$ equipped with the (idempotent) operation $x, y \mapsto
x \wedge y$ as addition.
\enddefinition

Note that the operations~$\vee = \oplus$ and~$\wedge$ (``union'', or
``addition'', and ``intersection'') are associative and commutative for
arbitrary sets of operands (summands or factors); compare, e.g.,~\cite{24,
Chap.~1, Sec.~6}.

\chead
2.4
\endchead
Let~$S$ be an arbitrary partially ordered set.
 This set is called
(ordinally) {\it complete\/} if any of its subsets, including the empty
one, has a greatest lower and a least upper bound.
 A set~$S$ is called {\it
boundedly complete\/} or {\it conditionally complete\/} if each of its
nonempty subsets that is bounded from above (from below) has a least upper
bound (respectively, a greatest lower bound).

\proclaim{Lemma 2.1}
If any subset~$S_1$ in~$S$ that is bounded from above has a least upper
bound~$\sup S_1$, then any subset~$S_2$ in~$S$ that is bounded from below
has a greatest lower bound~$\inf S_2$\rom; if any subset of~$S$ that is
bounded from below has a greatest lower bound in~$S$, then any subset
of~$S$ that is bounded from above has a least upper bound in~$S$.
\endproclaim

This lemma is well known (see, e.g.,~\cite{23, 24}).
 Note that
$$
\align
\inf S_2
&=\sup\{x\mid x\preccurlyeq s \text{ for all } s\in S_2\},
\tag2.1
\\
\sup S_1
&=\inf\{x\mid s\preccurlyeq x \text{ for all } s\in S_1\}.
\tag2.1$'$
\endalign
$$

Recall that a {\it cut\/} in~$S$ is a subset of the form $I(X)=\{s\in S\mid x\preccurlyeq s$
for all $x\in X\subset S\}$,
where~$X$ is any subset of~$S$.
 By~$\wh S$ denote the set of all cuts
endowed with the following partial order:
$I_1\preccurlyeq I_2$ if and only if
$I_1\supset I_2$; the set~$\wh S$ is called the {\it normal completion\/} of the ordered
set~$S$, see~\cite{23}.

The map $i\:x\mapsto I(\{x\})$
is an embedding of~$S$ into~$\wh S$ with the following well-known properties:

1. The embedding $i\:S\to \wh S$
preserves all greatest
lower and least upper bounds that exist in~$S$ (thus we can identify~$S$
with a subset of~$\wh S$).

2. The cut~$I(X)$ coincides with~$\sup(i(X))$ for any subset~$X$ of~$S$;
in particular, any element of the completion~$\wh S$
is the least upper
bound of some subset~$X$ of~$S$.

3.
If $S$ is complete, then $\wh S = S$; in particular,
$\wh{\!\wh S}= \wh S$
(here we take into account the
natural identification of~$\wh S$ with a subset of~$\wh{\!\wh S}$, see property~1).

4.
If $\{X_\alpha\}_{\alpha \in A}$ is a family of subsets of~$S$, then
$$
\sup_\alpha I(X_\alpha)
=I\biggl(\bigcup_\alpha X_\alpha\biggr)
$$
in~$\wh S$; note that $I(\cup X_\alpha) = \cap I(X_\alpha)$.

5.
The normal completion~$\wh S$ has the structure of idempotent
semigroup ($x \oplus y = \sup \{x, y\}$); if~$S$ is an idempotent
semigroup, then the embedding $S\to \wh S$
is a semigroup homomorphism.
 In this case we call~$\wh S$ the $a$-{\it completion\/}
of the idempotent semigroup~$S$.

The set~$\wh S_b$ of all cuts of the form~$I(X)$, where~$X$ runs over
all subsets of~$S$ that are bounded from above, is a conditionally complete
set; the set~$\wh S_b$ is a subsemigroup of~$\wh S$.
 If~$S$ is an
idempotent semigroup, then we call~$\widehat S_b$ the $b$-{\it completion\/}
of this idempotent semigroup.

It is clear that~$\widehat S_b$ consists of all elements of $\widehat S$
that are majorized by elements of~$S$; in particular, the cut
$I(\emptyset)\in \wh S_b$
(where~$\emptyset$ is the empty set)
is the zero of~$\wh S_b$.

\definition{Definition~2.4}
An idempotent semigroup~$S$ is called $a$-{\it complete\/} (or {\it
algebraically complete\/}) if it is complete as an ordered set.
 It is
called $b$-{\it complete\/} (or {\it boundedly algebraically complete\/})
if it is boundedly complete as an ordered set and has the neutral element
$\0$.
\enddefinition

It is easy to show that $\wh S=\wh S_b\cup\{\Bbb I_{\wh S}\}$;
note that $\Bbb I_{\wh S}=\sup\wh S$
may belong to~$S$.
 If~$S$
is a~$b$-complete idempotent semigroup, then~$\wh S_b$ coincides with~$S$.

\chead
2.5.
 Important notation
\endchead
 Let~$S$ be an idempotent semigroup.
 By~$\oplus X$ and~$\wedge X$,
respectively, denote~$\sup(X)$ and~$\inf(X)$ for any subset~$X$ of~$S$, if
these bounds exist, i.e., belong to~$S$.
 Thus in an~$a$-complete idempotent
semigroup the sum is defined for any subset; in a~$b$-complete idempotent
semigroup the sum is defined for any subset that is bounded from above,
including the empty one.

Let~$X$ be a subset of an idempotent semigroup~$S$.
 We introduce the
following notation:
$$
\align
\Up(X)
&=I(X)=\{y\in S\mid x\preccurlyeq y \text{ for all } x\in X\},
\\
\Low(X)
&=\{y\in S\mid y\preccurlyeq x \text{ for all }x\in X\}.
\endalign
$$

\chead
2.6
\endchead
 Let~$S$ and~$T$ be idempotent semigroups.

\definition{Definition~2.5-\rm{a}}
 Suppose the semigroups~$S$ and~$T$ are~$a$-complete.
 We call a homomorphism $g\:S\to T$
{\it algebraically continuous}, or an $a$-{\it homomorphism} for short, if
$$
g(\oplus X)=\oplus g(X),
\tag2.2
$$
i.e.,
$$
g\biggl(\sup_{x\in X}x\biggr)=\sup_{x\in X}g(X)
\tag2.2$'$
$$
for any subset~$X\subset S$.

For arbitrary idempotent semigroups $S$ and $T$ we call a homomorphism $g\:S\to T$
an~$a$-{\it homomorphism}, if it is uniquely extended to~$\wh S$ by an~$a$-homomorphism $\wh g\:\wh S\to \wh T$
of the corresponding normal completions.
\enddefinition

Many authors have introduced conditions similar to~\thetag{2.2$'$};
in~\cite{44}, this
condition is called {\it monotonic continuity}.

It is natural to consider the following variant of this definition
(in the spirit of~\cite{14}):

\definition{Definition~2.5-\rm{b}}
 Let~$S$ and~$T$ be~$b$-complete idempotent semigroups.
 A homomorphism
$g\:S\to T$ is called {\it boundedly algebraically continuous}, or a~$b$-{\it homomorphism} for short, if
condition~\thetag{2.2}, which coincides with~\thetag{2.2$'$},
is satisfied for any subset~$X \subset S$ bounded
from above.

For arbitrary idempotent semigroups~$S$ and~$T$, a homomorphism $g:S\to T$
is called a $b$-{\it homomorphism}, if it is uniquely extended to $\wh S_b$
by a~$b$-homomorphism
$\wh g\:\wh S_b\to \wh T_b$
of the corresponding~$b$-completions.
\enddefinition

Applying equality~\thetag{2.2} to the empty set, we
see that~$a$-homomorphisms and~$b$-homomor\-phisms take zero to zero
provided a zero element exists.  Note that in the general case
an~$a$-homomorphism $S \to T$ does not necessarily take the element $\Bbb
I_S=\sup S$ to $\Bbb I_T=\sup T$, even if both of these elements exist.

\proclaim{Proposition~2.1}
The composition \rom(i.e., the product\rom) of~$a$-homomorphisms
\rom($b$-homomorphisms\rom) is an~$a$-homomorphism
\rom(respectively, a~$b$-homomorphism\rom).
\endproclaim

This statement follows immediately from the definitions.

In fact, we are dealing with three different categories of idempotent
semigroups, where the morphisms are homomorphisms,
$a$-ho\-mo\-mor\-ph\-isms, and $b$\/-ho\-mo\-mor\-ph\-isms, respectively.
There exist other interesting categories of idempotent semigroups.
  The
notions introduced so far can be extended to the case of idempotent
semirings (see below).

By $\Hom(S_1,S_2)$ denote the set of all homomorphisms of an idempotent
semigroup~$S_1$ to an idempotent semigroup~$S_2$;
by~$\Hom_a(S_1,S_2)$ denote the set of all~$a$-homomorphisms.
 Finally, by $\Hom_b(S_1,S_2)$
denote the set of all~$b$-homomorphisms of a~$b$-complete idempotent
semigroup~$S_1$ to a~$b$-complete idempotent semigroup~$S_2$ (though the notion of~$b$-homomorphism can be extended to arbitrary idempotent
semigroups using the notion of conditional (bounded) completion~\cite{23}).

\proclaim{Proposition~2.2}
All three sets~$\Hom(S_1, S_2)$, $\Hom_a(S_1, S_2)$, and~$\Hom_b(S_1, S_2)$
are idempotent semigroups with respect to the pointwise sum. If~$S_2$ is
an~$a$-complete \rom(a~$b$-complete\rom) semigroup, then~$\Hom_a(S_1, S_2)$
\rom(respectively,~$\Hom_b(S_1, S_2)$\rom) is an~$a$-complete
\rom(respectively, a~$b$-complete\rom) idempotent semigroup.
\endproclaim

This statement follows directly from the definitions.


\chead
2.7
\endchead
We borrow the following definition from L.~Fuchs' book~\cite{27}; the
notation was introduced in Sec.~2.5 above.

\definition{Definition~2.6}
Let~$S$ be an idempotent semigroup.
 By definition, put~$X^\preccurlyeq = \Low(\Up(X))$ for any subset~$X$ of~$S$;
we call~$X^\preccurlyeq$ an~$o$-{\it closure\/} of~$X$.
\enddefinition

Note that the~$o$-closure of the empty set may be nonempty.

\definition{Definition~2.7}
Let~$S$ and~$T$ be idempotent semigroups.
 A map $f\:S\to T$
is said to be~$a$-{\it regular\/} if
$f(X^\preccurlyeq)\subset f(X)^\preccurlyeq$
for any subset~$X$ of~$S$;
this map is said to be~$b$-{\it regular\/}
if
$f(X^\preccurlyeq)\subset f(X)^\preccurlyeq$
for any subset~$X$ of~$S$ that is bounded from above.
\enddefinition

It is easy to prove that~$b$-regular (and therefore~$a$-regular) maps
are homomorphisms of idempotent semigroups.
  For a stronger statement, see
Proposition~2.3 below.

For example, let~$\Bbb R$ be the set of all real numbers equipped with the
structure of idempotent semigroup with respect to the operation $\oplus =
\max$.
 A map $f\:\Bbb R\to \Bbb R$ is a homomorphism of this
idempotent semigroup to itself if it is nondecreasing.
 This homomorphism is~$b$-regular (and~$a$-regular) if and only if it is lower semicontinuous
(for the definition of semicontinuity and some generalizations, see
Sec.~2.8 and Sec.~2.9 below).

The following statements can easily be checked.

\proclaim{Proposition~2.3}
Let~$S$ and~$T$ be idempotent semigroups.  A map $f\:S\to T$ is~$a$-regular
\rom($b$-regular\rom) and takes zero to zero if and only if it is
an~$a$-homomorphism \rom(respectively, a~$b$-homomorphism\rom).
\endproclaim

\proclaim{Proposition~2.4}
 A homomorphism~$f$ of an idempotent semigroup~$S$ to an idempotent
semigroup~$T$ is an~$a$-homomorphism if and only if it is
a~$b$-homomorphism and $\Up(f(X)) = \Up(f(S))$ for any subset~$X$ of~$S$
that is not bounded from above.
\endproclaim

\definition{Definition~2.8}
Let~$S,T$, and~$U$ be idempotent semigroups.
 We say that a map $f\:S\times T\to U$ is a {\it
 separate~$a$-homomorphism\/} ($b$-{\it homomorphism\/}) if the maps
$f(s,\,\cdot\,)\: t\mapsto f(s,t)$
and
$f(\,\cdot\,,t)\:s\mapsto f(s,t)$
are~$a$-homomorphisms
(respectively,~$b$-homomorphisms).
\enddefinition

In~\cite{45}, maps of this type were called bimorphisms.

\proclaim{Proposition 2.5}
Let $S,T$, and~$U$ be idempotent semigroups.

\rom{1}.
The normal completion $\wh{S \times T}$ of the direct product~$S \times T$
of the semigroups~$S$ and~$T$ coincides with the direct product~$\wh S
\times \wh T$ of their normal completions \rom(i.e., there exists a canonical
isomorphism that identifies these semigroups\rom).  An analogous statement
holds for $b$-completions.

\rom{2}.
The following equalities \rom(isomorphisms\rom) hold:
$$
\align
\oplus(X\times Y)
&=(\oplus X)\times(\oplus Y),
\tag2.3
\\
\Low(X\times Y)
&=\Low(X)\times\Low(Y),
\tag2.4
\\
\Up(X\times Y)
&=\Up(X)\times\Up(Y),
\tag2.5
\endalign
$$
where $X\subset S$, $Y\subset T$.

\rom{3}.
A map $f\:S\times T\to U$ is a separate~$a$-homomorphism
\rom($b$-homomorphism\rom) if and only if it can be \rom(uniquely\rom)
extended to~$\wh S$ \rom(respectively, to~$\wh S_b$\rom) by a
separate~$a$-homomorphism $\wh f\:\wh S\times\wh T\to \wh U$
\rom(respectively, a separate~$b$-homomorphism $\wh f\:\wh S_b\times\wh
T_b\to \wh U_b$\rom).
\endproclaim

Let us prove, for instance, Proposition~2.3 for the case in which~$S$
and~$T$ are~$a$-complete and $a$-regularity is considered.  Note that in
 an~$a$-complete idempotent semigroup we have
$$
X^\preccurlyeq=\Low(\oplus X)
\tag2.6
$$
for any subset~$X$.
 Let~$f$ be an~$a$-regular map.
 Using
equality~\thetag{2.6}, we get $f(\oplus X)\in f(X^\preccurlyeq)$; thus
$f(\oplus X)\in f(X^\preccurlyeq)\subset f(X)^\preccurlyeq
=\Low(\oplus f(X))$
and therefore
$f(\oplus X)\preccurlyeq\oplus f(X)$.
 Since any homomorphism preserves the order, we obtain
$f(\oplus X)\succcurlyeq\oplus f(X)$.
 Thus
$f(\oplus X)=\oplus f(X)$,
i.e.,~$f$ is an~$a$-homomorphism.
 Now let~$f$ be an~$a$-homomorphism; then
$f(\oplus X)=\oplus f(X)$.
Since~$f$ preserves the order, we can use formula~\thetag{2.6} to get
$$
f(X^\preccurlyeq)
=f(\Low(\oplus X))\subset\Low(f(\oplus X))
=\Low(\oplus f(X))=f(X)^\preccurlyeq,
$$
as claimed. The other statements can be proved similarly.

\chead
2.8.
 Semicontinuity
\endchead
The notion of semicontinuity is involved in the relationship between the
``topological'' approach, which  will be discussed in forthcoming
papers, and the ``algebraic'' one; also, this notion is needed here for the
presentation of some important examples.

By $\ol{\Bbb R}$ denote the extended real axis $\Bbb R\cup\{-\infty,+
\infty\}$ endowed with the usual order~$\le$ (which coincides with
the standard order in $\wh{\Bbb R}_{\max}$, so $\ol{\Bbb R}$ may be
identified with $\wh{\Bbb R}_{\max}$).  The set~$\ol{\Bbb R}$ has the
standard topology, which coincides with the ordinary topology on~$\Bbb
R$, so that the space~$\ol{\Bbb R}$ is homeomorphic to the segment $[-1, 1]$
of the real axis.

A function $f\:T\to\ol{\Bbb R}$ defined on an arbitrary
topological space~$T$ is said to be {\it lower semicontinuous\/}
\cite{46, Chap.~IV} if for any finite number~$s$ the set
$$
T_s=\{t\in T\mid f(t)\le s\}
$$
is closed in~$T$.
 This condition is equivalent to the condition that the
set $\{t\in T\mid f(t)> s\}$ is open in~$T$.
 An {\it upper
semicontinuous function\/} is defined similarly; a function~$f$ is upper
semicontinuous if and only if the function~$-f$ is lower semicontinuous.
Obviously, a real-valued function~$f$ is continuous if and only if it is
both upper and lower semicontinuous.
 The set
$$
\epi(f)=\{(t,s)\mid f(t)\le s\}
$$
in $T\times\ol{\Bbb R}$ is called the {\it epigraph\/} of the function~$f$.
 A function~$f$ taking finite values is lower semicontinuous if and
only if its epigraph is closed, see, e.g.,~\cite{47, 48}.
 Thus a subset in~$T$ is open (closed) if and only if its characteristic function is lower
(respectively, upper) semicontinuous.

Suppose $\{f_\alpha\}_{\alpha \in A}$ is a family of functions defined
on~$T$ and taking values in~$\ol{\Bbb R}$.  The {\it lower envelope\/}
$\inf_\alpha f_\alpha$ (the {\it upper envelope\/} $\sup_\alpha f_\alpha$)
is the function that is defined on $T$ and assigns the value
$\inf_\alpha(f_\alpha(t))$ (respectively, $\sup_{\alpha} (f_{\alpha} (t))$)
to each $t \in T$, see~\cite{46}.  It is well known that the {\it upper
envelope of a family\/}~$\{f_\alpha\}$ {\it of lower semicontinuous
functions on\/}~$T$ {\it is lower semicontinuous on}~$T$.  An analogous
statement holds for the lower envelope of a family of upper semicontinuous
functions, see~\cite{46, 47}.

The definition of semicontinuity can be extended to the case of maps
$f\: T\to S$, where~$T$ is a topological space and~$S$ is a (partially) ordered set.
 We call this map {\it lower semicontinuous\/}
if the set
$$
T_s=\{t\in T\mid f(t)\preccurlyeq s\}
\tag2.7
$$
is closed in~$T$ for all $s\in S$.
 A map is {\it upper
semicontinuous\/} if it is lower semicontinuous with respect to the dual
(opposite) order in~$S$.

 If~$S$ is a complete lattice (or, equivalently, an~$a$-complete idempotent
semigroup), then it can be equipped with the order topology in the sense
of~\cite{23, Chap.~10}.
 In~\cite{7, 8}, the notion of semicontinuity is
studied under some additional conditions on this topology.
 In~\cite{46, Chap.~IV, Sec.~6}, semicontinuous maps were studied for the
 case in which the order in~$S$ is linear (in particular, see exercise~6,
where~$S$ is compact and open intervals form a base for the topology
 in~$S$).  In these cases the classical properties of semicontinuity are
fulfilled.

\definition{Definition~2.9}
A {\it topological idempotent semigroup\/}~$S$ is an idempotent semigroup~$S$
endowed with a topology such that its subsemigroup $S^b=\{s\in S\mid s\preccurlyeq b\}$
is closed for any $b \in S$.
\enddefinition

\proclaim{Proposition~2.6-\rm{a}}
 Suppose~$T$ and~$S$ are~$b$-complete topological idempotent semigroups,~$S$ is~$a$-complete, and for any nonempty bounded subsemigroup~$X$ in~$T$ the
element~$\oplus X$ lies in the closure~$\ol X$ of the subset~$X$ in~$T$.
 Then a homomorphism $f\:T\to S$ taking~$\0_T$ to~$\0_S$ is an~$a$-homomorphism if and only if
the map~$f$ is lower semicontinuous.
\endproclaim

\demo{Proof}
First suppose that the map~$f$ is lower semicontinuous.
 For an arbitrary
bounded subsemigroup $X \subset T$ let $s = \oplus f(X)$.
 It is clear that
$X \subset T_s$, where~$T_s$ is defined by~\thetag{2.7}.
 Now it follows from the
fact that~$T_s$ is closed and the assumption of Proposition~2.6-a that
$\oplus X \in T_s$, i.e., $f(\oplus X)\preccurlyeq s=\oplus f(X)$.
  Since~$f$
preserves the order, the inequality $f(\oplus X)\succcurlyeq\oplus f(X)$ also holds, so $f(\oplus X)=\oplus f(X)$ if~$X$ is a nonempty bounded
subsemigroup in~$T$.
  Suppose now that~$X$ is an arbitrary nonempty subset
of~$T$ and denote by~$X^{\oplus}$ the subsemigroup of~$T$ generated by~$X$
(i.e., consisting of all finite sums of elements of~$X$).
 Since the map~$f$ is a homomorphism, $f(X^{\oplus}) = f(X)^{\oplus}$; taking into account
that $\oplus X = \oplus X^{\oplus}$, we obtain
$$
f(\oplus X) = f(\oplus X^{\oplus}) = \oplus f(X^{\oplus})
= \oplus (f(X)^{\oplus}) = \oplus f(X).
$$
It remains to prove that $f(\oplus\emptyset)=\oplus f(\emptyset)$,
but this means that~$\0_T$ is taken to~$\0_S$.
 This completes the proof of the
first part of Proposition~2.6-a.

Now suppose that the map~$f$ is an~$a$-homomorphism.
 We must prove that~$T_s$ is closed in~$T$ for any $s \in S$, i.e., that the limit of any net
(or generalized sequence) consisting of elements of~$T_s$ belongs to~$T_s$.
  Suppose that $x_\alpha \in T_s$ and $x = \lim_\alpha x_\alpha$; let
$t=\bigoplus_\alpha x_\alpha=\sup_\alpha x_\alpha$.
 Then
$f(t)=\bigoplus_\alpha f(x_\alpha)=\sup_\alpha f(x_\alpha)$.
 Hence $f(t)\preccurlyeq s$
(since $f(x_\alpha)\preccurlyeq s$),
$f(x)\preccurlyeq f(t)$
(since $x\preccurlyeq t$ and the homomorphism preserves the order),
and $x\in T_{f(t)}\subset T_s$, i.e., $x\in T_s$.
 Thus the homomorphism~$f$ is lower semicontinuous.
 This completes the proof.
\qed\enddemo

\proclaim{Proposition~2.6-\rm{b}}
 Suppose $f\:T\to S$ is a homomorphism of~$b$-complete topological
idempotent semigroups, $f(\0_T) = \0_S$, and for any nonempty subsemigroup~$X$ in~$T$ that is bounded from above the element~$\oplus X$ lies in the
closure of the set~$X$ in~$T$.
 Then~$f$ is a~$b$-homomorphism if and only
if its restriction to the closed subsemigroup
$T^x=\{t\in T\mid t\preccurlyeq x\}$ is lower semicontinuous for any~$x \in T$.
\endproclaim

Since~$T^x$ is an~$a$-complete idempotent semigroup and~$f(T^x)$ belongs to
the~$a$-complete idempotent semigroup $S^{f(x)}=\{s\in S\mid s\preccurlyeq f(x)\}$, Proposition~2.6-b is a straightforward consequence of
Proposition~2.6-a.
 The conditions of Proposition~2.6-b are fulfilled for
a wide class of idempotent semigroups; among the examples are all Banach
$L$-lattices (see~\cite{25, 26}), including the lattice of integrable
functions, spaces of semicontinuous real-valued functions with the `weak'
topology in the sense of~\cite{2, 6--8}, etc.

The analogy between the lower semicontinuity property for semirings of
the~$\Bbb R_{\max}$ type and continuity of countably additive measures in~$\emptyset$ was observed by A.~M.~Chebotarev in~\cite{5, Sec.~1.1.6} and applied to
the construction of the idempotent analog of the integral~\thetag{1.1} and
the ``Fourier--Legendre transform''~\thetag{1.5}.

\remark{Remark 2.2}
Suppose a topological idempotent semigroup~$T$ satisfies the conditions of
Proposition~2.6-a (2.6-b); then any of its subsemigroups also satisfies
these conditions whenever it is closed in the topology in~$T$
and is closed with respect to summation of arbitrary (respectively,
arbitrary bounded from above) subsets.
\endremark

\chead
2.9.
 Examples
\endchead
\dhead
\rom{2.9.1}
\enddhead
The set~$\Bbb R$ of real numbers is an idempotent semigroup with respect to
the operation ${x \oplus y} = \max \{x, y\} = \sup \{x, y\}$; here, the
standard order~$\preccurlyeq$
coincides with the usual order~$\le$.
 By~$\Bbb R_{\max}$ denote the idempotent monoid $\Bbb R\cup\{-\infty\}$ obtained by the
addition of the ``minimal'' element~$-\infty$ to~$\Bbb R$,
i.e., $\0= -\infty$.
The semigroup~$\Bbb R_{\max}$
is nothing but the~$b$-completion of the idempotent
semigroup~$\Bbb R$,
so~$\Bbb R_{\max}$ is~$b$-complete.
 The normal completion
$\wh{\Bbb R}_{\max}
=\Bbb R_{\max}\cup\{+\infty\}
=\Bbb R\cup\{-\infty\}\cup\{+ \infty\}$
is an~$a$-complete idempotent semigroup.
 The dual order~$(\ge)$ corresponds to the~$o$-dual semigroups with the operation $x
\oplus y = \min\{x, y\}$, which include the monoid $\Bbb R_{\min}=\Bbb
R\cup\{+ \infty\}$ and its normal completion $\wh{\Bbb R}_{\min} =\Bbb
R\cup\{+\infty\}\cup\{-\infty\}$; obviously, the semigroups~$\Bbb R_{\max}$
and~$\Bbb R_{\min}$ are isomorphic, as well as their completions.

In the sequel, as in other examples (see below), $\Bbb R_{\max},\Bbb R_{\min}$, and
their completions will have the structure of idempotent semirings in
addition to that of idempotent semigroup.

\dhead
\rom{2.9.2}
\enddhead
The completion $\wh{\Bbb R}_{\max}$ is isomorphic to the semigroup formed by the
segment $[0, 1]$ under the operation $x \oplus y = \max \{x, y\}$.

\dhead
\rom{2.9.3}
\enddhead
 The set~$\Bbb Z$ of all integers is an idempotent semigroup with respect to
the operation $x \oplus y = \max \{x, y\}$.
 By $\Bbb Z_{\max}$ denote a~$b$-complete idempotent semigroup $\Bbb Z\cup\{-\infty\}$ obtained by
adding the element $\0 = -\infty$.
 After adding the largest element $\Bbb I=+ \infty$ to
$\Bbb Z_{\max}$, we get its normal completion
$\wh{\Bbb Z}_{\max}$.
These idempotent semigroups are subsemigroups of the idempotent semigroups
$\Bbb R,\Bbb R_{\max},\wh{\Bbb R}_{\max}$, of Example~2.9.1. Similarly are
defined the idempotent semigroups~$\Bbb Z_{\min}$ and $\wh{\Bbb Z}_{\min}$.

\dhead
\rom{2.9.4}
\enddhead
Let~$X$ be an arbitrary set, $S$~be an idempotent semigroup.
 By~$\Map(X,S)$ denote the set of all functions defined on~$X$ and taking values in~$S$
(i.e., the set of all maps~$X \to S$).
 The standard order on the set~$\Map(X, S)$ is defined by the rule
$$
f\preccurlyeq g,\ \text{if and only if}\
f(x)\preccurlyeq g(x)\ \text{for all}\ x\in X.
\tag2.8
$$
This order corresponds to the operation of pointwise addition $f,g\mapsto
f\oplus g$ that endows $\Map(X,S)$ with the structure of idempotent
semigroup.

\dhead
\rom{2.9.5}
\enddhead
 By~$\sB(X,S)$ denote the subsemigroup in~$\Map(X,S)$ formed by
{\it bounded\/} functions (i.e., functions with bounded range).
 This
idempotent semigroup is~$b$-complete if~$S$ is~$b$-complete.

\dhead
\rom{2.9.6}
\enddhead
 By~$\sB(X)$ denote the semigroup~$\sB(X,S)$ defined in Example~2.9.5,
if~$S$ coincides with the idempotent semigroup~$\Bbb R$ defined in
Example~2.9.1.

\dhead
\rom{2.9.7}
\enddhead
By~$C(X)$ denote the set of all continuous real-valued functions on a
topological space~$X$.
 The operation $f \oplus g$ is defined by the
equality
$$
(f\oplus g)(x)=\max\{f(x),g(x)\}
\tag2.9
$$
for all $x\in X$;
this operation turns~$C(X)$ to an idempotent semigroup.

\dhead
\rom{2.9.8}
\enddhead
By~$USC(X)$ ($LSC(X)$) denote the set of upper (respectively,
lower) semicontinuous real-valued functions on a topological space~$X$.
Recall that a function~$f(x)$ is upper semicontinuous if
$$
f(x)=\inf_\nu\{f_\nu(x)\}
\tag2.10
$$
for all $x\in X$ (the pointwise~$\inf$), where
$\{f_\nu\}$ is any family
of upper semicontinuous (in particular, continuous) functions (see
Sec.~2.8 above).
 Similarly,~$f(x)$ is lower semicontinuous if
$$
f(x)=\sup_\nu\{f_\nu(x)\}
\tag2.10$'$
$$
for all $x\in X$, where
$\{f_\nu\}$ is a family of lower semicontinuous
functions.

Let~$\Bbb R$ be the idempotent semigroup of real numbers with respect to
the operation~$\oplus = \max$ (Example~2.9.1).
 Then $USC(X)$ and $LSC(X)$ are
subsets of the partially ordered set $\Map(X, \Bbb R)$.
 Thus they are
partially ordered by the induced order
$$
f\le g\iff f(x)\le g(x)
\tag2.11
$$
for all $x\in X$.
 The idempotent operation
$$
f\oplus g=\sup\{f,g\}
\tag2.12
$$
endows both $USC(X)$ and $LSC(X)$ with the structure of idempotent
semigroups (subsemigroups of $\Map(X, \Bbb R)$).

\remark{Remark~2.3}
In general, least upper (greatest lower) bounds of infinite sets in the
semigroup $\Map(X, \Bbb R)$ and its subsemigroups $USC(X)$ and $LSC(X)$ do
not coincide.
 In $LSC(X)$, only least upper bounds coincide with those in
$\Map(X, \Bbb R)$; in $USC(X)$, only greatest lower bounds do.
 Any set in
$USC(X)$ that is bounded from below has a greatest lower bound (the
pointwise~$\inf$); thus it follows from Lemma~2.1 that the idempotent
semigroup $USC(X)$ is a~$b$-complete lattice.
 By the same argument,
$LSC(X)$ is a~$b$-complete lattice.
\endremark

\dhead
\rom{2.9.9}
\enddhead
The set $\Conv(X, \Bbb R)$ of all convex real-valued functions defined on a
convex subset~$X$ of some vector space\footnote{Recall that a function is
said to be {\it convex\/} if its epigraph (for the definition see
Sec.~2.8 above) is a convex subset of $X \times \Bbb R$; a function~$f$
is said to be {\it concave\/} if the function~$-f$ is convex (see~\cite{47,
48}).} is an idempotent semigroup (and a~$b$-complete lattice) with
respect to the operation~\thetag{2.12}.
 This idempotent semigroup is a
subsemigroup of $\Map(X, \Bbb R)$, and the embedding is an~$a$-homomorphism.

\dhead
\rom{2.9.10}
\enddhead
The set $\Conc(X, \Bbb R)$ of all concave real-valued functions defined on a
convex subset~$X$ of some vector space has the structure of idempotent
semigroup under the restriction of the order defined by~\thetag{2.11}.
This idempotent semigroup is a lattice but is not a subsemigroup of
$\Map(X, \Bbb R)$.

\dhead
\rom{2.9.11}
\enddhead
Let~$X$ be a measure space with measure~$\mu$ and $L^1(X, \mu)$ be the
(Banach) space of~$\mu$-integrable functions on this space (i.e., classes
of functions on~$X$ that differ on a set of zero measure~$\mu$).
 Then
$L^1(X, \mu)$ is a boundedly complete idempotent semigroup with respect to
the operation~\thetag{2.12}
generated by the standard order~\thetag{2.11}, where the
inequality $f(x)\le g(x)$ is supposed to hold almost everywhere.
This result is true for any Banach space $L^p(X, \mu)$ with $1 \le p
\le \infty$.
 These spaces provide an example of boundedly complete
Banach lattices (for the theory of vector lattices see~\cite{25, 26} and
\cite{23, 24}).

\dhead
\rom{2.9.12}
\enddhead
Suppose~$X$ is a metric space with metric~$\rho$; by $\lip(X)$ denote the
semigroup of all real-valued functions on~$X$ satisfying the {\it Lipschitz
condition}
$$
|f(x)-f(y)|\le\rho(x,y).
\tag2.13
$$
The structure of idempotent semigroup is defined by~\thetag{2.11}
and~\thetag{2.12}.

By $\blip(X)$ denote this semigroup supplemented with the element~$\0$,
so that $\blip(X)$ is an idempotent monoid.

\remark{Remark~2.4}
The semigroups $USC(X),LSC(X),\Conv(X,\Bbb R),\lip(X)$, and $L^p(X,\mu)$
are conditionally complete (but neither~$b$-complete nor~$a$-complete).
Adding~$\0$ to any of these semigroups, we obtain $b$-complete
semigroups.
\endremark

\ahead
3.
 IDEMPOTENT SEMIRINGS, SEMIFIELDS, AND QUASIFIELDS
\endahead

\bhead
3.1.
 Idempotent semirings and semifields
\endbhead
Idempotent semirings and semifields are the main objects of Idempotent
Mathematics.

\definition{Definition~3.1}
An {\it idempotent semiring\/} (or {\it semiring\/} for short) is an
idempotent semigroup~$K$ (with addition~$\oplus$) endowed with an
additional associative multiplication operation~$\odot$ such that for all
$x, y, z \in K$
$$
\align
x\odot(y\oplus z)
&=(x\odot y)\oplus(x\odot z),
\tag3.1
\\
(y\oplus z)\odot x
&=(y\odot x)\oplus(z\odot x).
\tag3.1$'$
\endalign
$$
\enddefinition

Note that $a\odot x\preccurlyeq a\odot y$ if
$x\preccurlyeq y$ for all $x,y,a\in K$.

An element $\1 \in K$ is {\it unit\/} of a semiring $K$ if it is
neutral with respect to multiplication, i.e., if
$$
\1\odot x=x\odot\1=x
\tag3.2
$$
for all $x\in K$.
 In the sequel {\it we always assume that all idempotent
semirings contain a unit}, unless otherwise stated.

An element $\0 \in K$ is {\it zero\/} of a semiring~$K$ if it is zero
with respect to addition~$\oplus$, i.e., $x \oplus \0 = x$, and
$$
x\odot\0=\0\odot x=\0
\tag3.3
$$
for all $x\in K$.
 We also denote zero of a semiring~$K$ by~$\0_K$.
 An
idempotent semiring with zero is often called a {\it dioid}, see,
e.g.,~\cite{12}.
 We always assume that $\0\ne \1$, unless otherwise stated.
  If a
zero exists in a semiring, it is also a zero of the additive semigroup of
this semiring.
 The converse is not always true.
 The set $\{1,2,3,\dots\}$ of natural numbers is a semiring with respect to the operations
$\oplus = \max$, $\odot = +$; the element~$1$ is a zero of its additive
semigroup but not of the whole semiring.

An idempotent semiring~$K$ is {\it commutative\/} if the multiplication
operation is commutative.
 There exist different versions of the axiomatics of
idempotent semirings; for these, as well as for historical remarks,
see~\cite{8, 9, 12, 13, 16--21, 49}.

\definition{Definition~3.2}
An {\it idempotent division semiring\/} is an idempotent semiring with
unit in which any nonzero element has a multiplicative
inverse.
 An {\it idempotent semifield\/} (or {\it semifield\/} for short)
is a commutative idempotent division semiring.
\enddefinition

This notion of idempotent semifield does not coincide with the notion of
semifield in the sense of~\cite{50}.

\bhead
3.2. Examples
\endbhead
\dhead
\rom{3.2.1}
\enddhead
All examples listed in Sec.~2.9 above, except Example~2.9.12, become
idempotent semirings after a suitable product is defined; a semiring
without unit appears in Example~2.9.12 if~$\odot$ is defined as~$\min$.
 In
Example~2.9.1, multiplication~$\odot$ is defined as the ordinary
addition~$+$ and extended to the additional elements~$\pm \infty$ in the
natural way.  For example, in $\wh{\Bbb R}_{\max}$, equalities~\thetag{3.3}
and the rules $x\odot(+\infty)=x+(+\infty)=+\infty =\sup\wh{\Bbb R}_{\max}=
\Bbb I$ hold for all $x\ne \0$.  The element~$\1$ coincides with the
 usual zero.  The idempotent semiring~$\Bbb R$ with these operations
is often denoted by $\Bbb R(\max, +)$.  In Example~2.9.2, multiplication is
induced by the isomorphism with $\wh{\Bbb R}_{\max}$.

The semigroup $\Map(X, S)$ is a semiring if~$S$ is an idempotent semiring;
in this case multiplication of functions is defined pointwise:
$$
(f\odot g)(x)=f(x)\odot g(x)
\tag3.4
$$
for any $x\in X$.
 The multiplication~$\odot$ in Examples~2.9.4--2.9.11 is defined in a
similar way.
 Unless otherwise stated, the multiplication $\odot$ of
real-valued functions is defined as
$$
(f\odot g)(x)=f(x)+g(x)
\tag3.5
$$
for any $x\in X$.

\dhead
\rom{3.2.2}
\enddhead
Suppose $S$ is an arbitrary idempotent semigroup and $\Hom(S, S)$ is the
semigroup of homomorphisms of $S$ to itself (i.e., endomorphisms).
The operation of composition of maps endows $\Hom(S, S)$ with the
structure of an idempotent semiring.

\bhead
3.3.
 Complete semirings
\endbhead
\dhead
\rom{3.3.1}
\enddhead
We begin with the definition of algebraic completeness.

\definition{Definition~3.3}
An idempotent semiring~$K$ is~$a$-{\it complete\/} ({\it algebraically
complete\/}) if it is an~$a$-complete idempotent semigroup and for any subset $\{x_\alpha\}$ of~$K$ and any $y \in K$,
$$
\biggl(\bigoplus_\alpha x_\alpha\biggr)\odot y
=\bigoplus_\alpha(x_\alpha\odot y),\qquad
y\odot\biggl(\bigoplus_\alpha x_\alpha\biggr)
=\bigoplus_\alpha(y\odot x_\alpha),
\tag3.6
$$
i.e.,
$$
\biggl(\sup_\alpha\{x_\alpha\}\biggr)\odot y
=\sup_\alpha\{x_\alpha\odot y\},\qquad
y\odot\biggl(\sup_\alpha\{x_\alpha\}\biggr)
=\sup_\alpha\{y\odot x_\alpha\}.
\tag3.6$'$
$$
\enddefinition

This means that the ``homotheties'' $x \mapsto x \odot y$ and $x \mapsto y
\odot x$ are $a$-homomorphisms of the semigroup $K$.

A semiring $K$ with zero and unit is $a$-complete if and only if it is a
complete dioid in the sense of~\cite{12}; the notion of algebraically
complete semiring dates back to~\cite{51} and is closely related to
different versions of the notion of Kleene algebra (see~\cite{52}).
 The
notion of $a$-complete idempotent semiring is essentially equivalent to
that of {\it quantale\/} (see, e.g.,~\cite{53}), which appears in
connection with the foundations of Quantum Theory.
 Morphisms of the
corresponding category {\bf Quant} preserve infinite sums.
 A quantale does
not necessarily contain a (two-sided) unit element; a quantale with unit is
called a {\it uniquantale}.
 Commutative $a$-complete idempotent semirings
were studied in~\cite{54}, where they were called commutative monoids in
the category of complete sup-lattices.

Note that an $a$-complete semiring $K$ cannot be a semifield (except for the
case $K = \{\0,\1\}$)~\cite{49, 14}.
 For instance, the element
$+ \infty$ is not invertible in~$\wh{\Bbb R}_{\max}$.

\dhead
\rom{3.3.2}
\enddhead
A weaker but more useful version of the notion of~$a$-completeness is the
bounded completeness or~$b$-completeness of semirings.

\definition{Definition~3.4}
An idempotent semiring~$K$ is called~$b$-{\it complete\/} (or {\it boundedly
complete\/}) if it is a~$b$-complete idempotent semigroup and
equalities~\thetag{3.6} hold for any subset~$\{x_\alpha\}$ of~$K$ that is bounded
from above and for any $y \in K$.
\enddefinition

{\it For simplicity, we shall always assume in the sequel that any~$b$-complete idempotent semiring has a (semiring) zero element, unless
otherwise stated}.
 In~\cite{14} commutative idempotent semirings with this
property are called {\it regular\/}; in~\cite{49} they are called boundedly
complete (BC) dioids.

Obviously,~$a$-completeness implies~$b$-completeness.

\dhead
\rom{3.3.3}
\enddhead
We introduce $b$-complete semirings with the following important additional
property.

\definition{Definition~3.5}
We call a semiring~$K$ {\it lattice~$b$-complete\/} if it is a~$b$-complete
idempotent semiring with zero and for any nonempty subset~$\{x_\alpha\}$
of~$K$ and any $y\in K\setminus\Bbb I$
$$
(\wedge_\alpha x_\alpha)\odot y
=\wedge_\alpha(x_\alpha\odot y),\qquad
y\odot(\wedge_\alpha x_\alpha)
=\wedge_\alpha(y\odot x_\alpha).
\tag3.7
$$
\enddefinition

Note that any nonempty subset of~$K$ is bounded from below by zero and
therefore has a greatest lower bound (which is a least upper bound
of the set of its lower bounds, see Lemma~2.1).

\proclaim{Proposition~3.1}
Any idempotent division semiring that is~$b$-complete as a semigroup
\rom(hence, any~$b$-complete idempotent semifield\rom) is
lattice~$b$-complete.
\endproclaim

\demo{Proof}
It is sufficient to prove equalities~\thetag{3.6} and~\thetag{3.7}.
 These equalities hold
trivially if~$y = \0$.
If $y\ne \0$, then~$y$ is invertible; thus
the homotheties $x \mapsto x \odot y$, $x \mapsto y \odot x$ are invertible
maps that preserve the order and, in particular, least upper and greatest
lower bounds. This means that~\thetag{3.7} and \thetag{3.6} hold.
\qed\enddemo

\dhead
\rom{3.3.4}
\enddhead
Suppose~$S$ is a semigroup (not necessarily idempotent) with respect to the
multiplication $x, y \mapsto xy$.  This semigroup is called {\it lattice
ordered\/} (see~\cite{23}) if~$S$ is a lattice and the semigroup
translations $x \mapsto xy$, $x \mapsto yx$ are isotonic (i.e.,
order-preserving).  By~$S_{\0}$ denote either the set~$S$ itself,
if it contains the least element~$\0$ such that $\0 x = x \0 = \0$ for all
$x \in S$, or the set $S_{\0}=S\cup\{\0\}$, where~$\0$ is an additional
element, in the converse case.  Let $x \oplus y = \sup\{x, y\}$ if $x, y
\in S$; also, let $x \oplus \0 = \0 \oplus x = x$ for all $x \in S_{\0}$.

\proclaim{Proposition~3.2}
The set~$S_{\0}$ is an idempotent semiring with respect to the addition $x
\oplus y$ and the multiplication $x \odot y = xy$, where by definition $x
\odot \0 = \0 \odot x = \0$ for all $x \in S_{\0}$.
\endproclaim

\proclaim{Proposition~3.3}
If~$G$ is a boundedly complete lattice ordered
group, then~$G_{\0}$ is a~$b$-complete idempotent division semiring, i.e.,
any nonzero element of the semiring~$G_{\0}$ has a multiplicative inverse.
Any~$b$-complete idempotent division semiring has the form~$G_{\0}$, where~$G$ is a boundedly complete lattice ordered group.
\endproclaim

These statements follow directly from the definitions.

\proclaim{Proposition~3.4}
Any~$b$-complete idempotent division semiring is commutative, i.e., is an
idempotent semifield.
\endproclaim

This statement follows from Proposition~3.3 and~\cite{27, Part~1, Chap.~5,
Theorem~18}.

\bhead
3.4.
 Quasifields
\endbhead
Let~$K$ be an idempotent semiring.

\definition{Definition~3.6}
An element $x \in K$ is {\it quasi-invertible\/} if there exists a set~$X$
of invertible elements of~$K$ such that $x = \oplus X$.

We call a semiring~$K$ a {\it division quasiring} if any of its nonzero
elements is quasi-invertible.
\enddefinition

\definition{Definition~3.7}
A semiring~$K$ is {\it integrally closed\/} if $x\preccurlyeq\1$ whenever the set $\{x^n\mid n=1,2,\dots\}$ of all powers of an element $x
\in K$ is bounded from above (for all $x \in K$).
\enddefinition

For division semirings, integral closedness means that the group of all
invertible elements is integrally closed in the sense of~\cite{23,
Chap.~XIII, Sec.~2} (sometimes this property is called complete integral
closedness).

\definition{Definition~3.8}
An idempotent division quasiring~$K$ is a {\it quasifield\/} if it is
integrally closed.
\enddefinition

The following statement follows directly from the definitions.

\proclaim{Proposition~3.5}
Any integrally closed semifield is a quasifield.
\endproclaim

Note that any quasifield is commutative, although this is not mentioned in
Definition~3.8 explicitly (see Proposition~3.7 below).

\bhead
3.5.
 Completion of semirings
\endbhead
Suppose~$K$ is an idempotent semiring.
 Considering it as an idempotent semigroup, we can construct its~$a$-completion~$\wh K$ and~$b$-completion~$\wh K_b$. In this section we discuss whether
it is possible to endow~$\wh K$ ($\wh K_b$) with the
structure of~$a$-complete (respectively,~$b$-complete) idempotent semiring.
By $\mu$ denote the map $\mu\:K\times K\to K$,
defined by the multiplication operation, i.e.,
$$
\mu\:(x,y)\mapsto x\odot y.
$$

\definition{Definition~3.9}
A semiring~$K$ is said to be~$a$-{\it regular\/} ($b$-{\it regular\/}) if
the map~$\mu$ has a unique extension $\wh\mu\:\wh K\times\wh K\to \wh K$
(respectively,
$\wh\mu\:\wh K_b\times\wh K_b\to \wh K_b$) that defines the
structure of~$a$-complete (respectively,~$b$-complete) idempotent
semiring in~$\wh K$ (respectively,~$\wh K_b$).
\enddefinition

Clearly, the~$a$-regularity ($b$-regularity) of~$K$ implies that the
multiplication~$\mu$ is a separate~$a$-homomorphism (respectively, a
separate~$b$-homomorphism) in the sense of Definition~2.8.

\definition{Definition~3.10}
If~$K$ is an~$a$-regular (a~$b$-regular) semiring, then the semiring~$\wh K$ (respectively,~$\wh K_b$) is called the~$a$-{\it
completion\/} (respectively,~{\it $b$-completion\/}) of the semiring~$K$;
we denote~$\widehat\mu(x, y)$ by $x \odot y$ as before.
\enddefinition

\proclaim{Proposition~3.6}
A semiring~$K$ with zero is~$a$-regular \rom($b$-regular\rom) if and only
if the homotheties
$$
\mu(\,\cdot\,,y)\:x\mapsto x\odot y,
\qquad
\mu(y,\,\cdot\,)\:x\mapsto y\odot x
\tag3.8
$$
are~$a$-homomorphisms \rom(respectively,~$b$-homomorphisms\rom)
of the additive semigroup of~$K$ for all $y\in K$.
\endproclaim

\demo{Proof}
If~$K$ is regular, then the homotheties~\thetag{3.8} are obviously regular.
 Suppose
now that the homotheties~\thetag{3.8} are~$a$-homomorphisms ($b$-homomorphisms); then by
Definition~2.5 they can be uniquely extended to the whole set~$\wh K$
(respectively,~$\wh K_b$).
 Thus the product $x \odot y$ is well
defined if one of the factors belongs to the semigroup completion of~$K$.
In the general case, let $x=\sup_\alpha x_\alpha$,
$y=\sup_\beta y_\beta$,
where $x_\alpha,y_\beta\in K$; let
$$
x\odot y
=\sup_\alpha(x_\alpha\odot y)
=\sup_\beta(x\odot y_\beta).
\tag3.9
$$
It is easily shown (by using Propositions~2.5 and 2.2) that the
product~\thetag{3.9} is well defined.  A straightforward calculation shows
 that this product is associative and endows~$\wh K$ (or~$\wh K_b$) with
the required structure.
\qed\enddemo

Recall that the~$a$-completion of a semiring~$K$ cannot be a semifield (if
$K\ne\{\0,\1\}$).
 Obviously, completions of commutative semirings are
commutative.

\proclaim{Proposition~3.7}
Any quasifield~$K$ is commutative and~$b$-regular\rom; hence~$\wh K_b$ is
a~$b$-complete idempotent semifield.
\endproclaim

\demo{Proof}
First we use Proposition~3.6 to check that the quasifield~$K$ is~$b$-regular.
  If an element $y \in K$ is invertible, then it follows from
Proposition~3.1 that the homotheties~\thetag{3.8} are~$b$-regular (i.e., are~$b$-homomorphisms).
 For a quasi-invertible $y \in K$ the same statement
follows from Proposition~2.2 (the corresponding homothety is the sum, or the
least upper bound, of a set of~$b$-homomorphisms).
 Since any nonzero
element of~$K$ is quasi-invertible, this implies that~$K$ is~$b$-regular
and~$\wh K_b$ is a~$b$-complete idempotent semiring.
 Now from the
results of~\cite{27, Chap.~V} (see Theorem~19 and its corollaries) it
follows immediately that~$\wh K_b$ is a~$b$-complete semifield.

It can easily be checked that any~$b$-complete semifield is integrally
closed and is a quasifield.
 Thus $b$-{\it complete semifields are the most
important examples of quasifields\/}; on the other hand, {\it any
quasifield~$K$ can be embedded in a~$b$-complete semifield\/}~$\wh K_b$.
\qed\enddemo

\remark{Remark~3.1}
If~$K$ is a~$b$-complete semifield that does not coincide with~$\{\0,\1\}$,
then $\wh K= K\cup\{\Bbb I\}$ (where $\Bbb I=\sup\wh K$, see Sec.~2.2
above) has the structure of commutative idempotent semiring, where $\0
\odot \Bbb I = \0$ and $x \odot \Bbb I = \Bbb I$ for all $x \ne \0$ in~$K$.  {\it Thus
any~$b$-regular quasifield is~$a$-regular}.
\endremark

\bhead
3.6.
 Examples
\endbhead
\dhead
\rom{3.6.1}
\enddhead
It is easy to check that {\it the semirings~$\Bbb R_{\max}$ and~$\Bbb R_{\min}$ are~$b$-complete idempotent semifields\rom;
the idempotent semiring~$\Bbb R$
is a
semifield and a quasifield\/} (see Examples~2.9.1 and 3.2.1).

\dhead
\rom{3.6.2}
\enddhead
The subsemiring in~$\Bbb R_{\max}$ formed by all integers and the element~$\0$ is a~$b$-complete semifield.

\dhead
\rom{3.6.3}
\enddhead
 The semiring
$K=\Bbb R\cup\{-\infty\}\cup\{+ \infty\}$
with the
operations $\oplus = \max$ and $\odot = \min$, where $\0 = -\infty$, $\1 =
+\infty$, is lattice~$b$-complete.

\dhead
\rom{3.6.4}
\enddhead
The semiring $\Map(X,\wh{\Bbb R}_{\max})$
of all functions taking values in
$\wh{\Bbb R}_{\max}$
(see Examples~2.9.3 and 3.2.1) is~$a$-complete but not lattice~$b$-complete.

\dhead
\rom{3.6.5}
\enddhead
It follows from Proposition~3.3 that any (boundedly) complete vector
lattice over~$\Bbb R$ generates a~$b$-complete idempotent semifield.
 In
particular, $L^p(x, \mu)_{\0} = L^p(x, \mu) \cup \{\0\}$ is a~$b$-complete semifield, see Examples~2.9.11, 3.2.1.

\dhead
\rom{3.6.6}
\enddhead
Any vector lattice can be considered as a semifield; in particular, the set~$C(X)$ of all continuous real-valued functions on a topological space~$X$
is a semifield and a quasifield with respect to the standard idempotent
operations defined in Examples~2.9.7 and 3.2.1.

\ahead
4.
 IDEMPOTENT SEMIMODULES AND SPACES
\endahead

\bhead
4.1.
 Basic definitions
\endbhead

\definition{Definition~4.1}
Let~$V$ be an idempotent semigroup, $K$ be an idempotent semiring, and let
a multiplication operation $k, x \mapsto k \odot x$ be defined so that the
following equalities hold
$$
\gather
(k_1\odot k_2)\odot x=k_1\odot(k_2\odot x),
\tag4.1
\\
(k_1\oplus k_2)\odot x=(k_1\odot x)\oplus(k_2\odot x),
\tag4.2
\\
k\odot(x\oplus y)=(k\odot x)\oplus(k\odot y),
\tag4.3
\\
\1\odot x=x
\tag4.4
\endgather
$$
for all $x,y\in V$, $k,k_1,k_2\in K$; then the semigroup~$V$ is called
a (left) {\it idempotent semimodule\/} over the semiring~$K$ (or simply a
{\it semimodule\/}).
\enddefinition

\definition{Definition~4.2}
Let~$V$ be an idempotent semimodule over~$K$.
 An element $\0_V \in V$ is
{\it zero\/} of the semimodule~$V$ if $k \odot \0_V = \0_V$ and $\0_V
\oplus x = x$ for all $k \in K$, $x \in V$.
\enddefinition

In the sequel, the zero of a semimodule~$V$ is denoted by~$\0$ if it exists
and this notation does not lead to confusion.

As usual, a {\it subsemimodule\/} of a semimodule~$V$ is a subsemigroup of~$V$ that is invariant under multiplication by coefficients from~$K$.

By $L_k$ denote the operator of ``homothety'' in $V$, i.e., $L_k\:x\mapsto k\odot x$, where $k\in K$, $x\in V$.
 Evidently, the map $k
\mapsto L_k$ is a homomorphism of~$K$ to $\Hom(V, V)$, where $\Hom(V, V)$
is the semiring of all homomorphisms of the semimodule~$V$ (see
Example~3.2.2 above).

\bhead
4.2.
 Complete and standard semimodules
\endbhead

\definition{Definition~4.3}
A semimodule~$V$ over an idempotent semiring~$K$ is called~$a$-{\it
complete\/} if~$V$ is an~$a$-complete idempotent semigroup, all homotheties~$L_k$ are~$a$-homomorphisms, the homomorphism $k \mapsto L_k$ is
(uniquely) extended to give a semigroup homomorphism~$\wh K \to
\Hom(V, V)$ defining a multiplication operation $\wh K \times V \to
V$, and the following equalities hold:
$$
(\oplus Q)\odot x=\oplus(Q\odot x),
\tag4.5
$$
i.e.,
$$
(\sup Q)\odot x=\sup_{k\in Q}\{k\odot x\}
\tag4.5$'$
$$
for all $x\in V$, $Q\subset \wh K$.
 The notion of~$b$-{\it complete\/} semimodule over~$K$ is defined similarly (but~$Q$ is
assumed to be a subset of~$\wh K_b$).
\enddefinition

\remark{Remark 4.1}
Assuming that the set~$Q$ is empty in~\thetag{4.5}, we see that
for all $x \in V$
$$
\0_K\odot x=\0_V.
\tag4.6
$$
\endremark

Note that idempotent semimodules were considered in many publications, see,
e.g.,~\cite{4--14, 16--21, 42--44,~49}.
 In~\cite{45},~$a$-complete
idempotent semimodules over commutative~$a$-complete idempotent semirings
were considered in different terms.
 For additional references see Sec.~6
below.

\definition{Definition~4.4}
A semimodule~$V$ over a lattice~$b$-complete idempotent semiring~$K$ (see
Definition~3.5) is called {\it standard\/} if~$V$ is~$b$-complete (in
particular,~$a$-complete) and for any $x \in V$ such that $x \ne \Bbb I = \sup V$ and any nonempty subset $Q \subset K$ the following equalities hold:
$$
(\wedge Q)\odot x=\wedge(Q\odot x),
\tag4.7
$$
i.e.,
$$
(\inf Q)\odot x=\inf_{k\in Q}\{k\odot x\}.
\tag4.7$'$
$$
\enddefinition

Note that~$Q$ is bounded from below since a~$b$-complete semiring~$K$ has
zero~$\0_K$.

\definition{Definition~4.5}
A semimodule~$V$ over an idempotent semiring~$K$ is called~$a$-{\it
regular\/} ($b$-{\it regular\/}) if the product $\mu\:K\times V\to V$ is a separate~$a$-homomorphism (respectively,~$b$-homomorphism)
of
idempotent semigroups in the sense of Definition~2.8.
\enddefinition

Using Proposition~2.5 and arguing as in Sec.~3.5, we obtain the
following statement.

\proclaim{Proposition~4.1}
Let~$V$ be an~$a$-regular \rom(a~$b$-regular\rom) semimodule over
a~$b$-regular semiring~$K$ and $\mu\:K\times V\to V$ be the corresponding
product. Then~$\mu$ has a unique extension $\wh\mu\:\wh K\times\wh V\to
\wh V$ \rom(respectively, $\wh\mu\:\wh K_b\times\wh V_b\to \wh V_b$\rom),
which defines on~$\wh V$ \rom(respectively, on~$\wh V_b$\rom)
the structure of~$a$-complete \rom(respectively,~$b$-complete\rom)
semimodule over~$K$ and~$\wh K_b$. If the semiring~$K$ is~$a$-regular, i.e.,
if~$\wh K$ is a semiring, then~$\wh V$ is an~$a$-complete semimodule
over~$\wh K$ as well.
\endproclaim

We denote $\wh\mu(k,x)$ by $k \odot x$ as before.

\definition{Definition~4.6}
If~$V$ is an~$a$-regular (a~$b$-regular) semimodule over a~$b$-regular
semiring~$K$, then we call the semimodule~$\wh V$ (respectively, $\wh V_b$)
over~$\wh K_b$ the~$a$-{\it completion\/} (respectively, the~$b$-{\it
completion\/}) of the semimodule~$V$ (see Proposition~4.1).
\enddefinition

We stress that the procedure of completion of a semimodule~$V$ requires that
the semiring~$K$ be replaced with its completion~$\wh K_b$.

\bhead
4.3.
 Idempotent spaces
\endbhead
To obtain substantial results on semimodules as well as maps and
functionals defined on semimodules, it is appropriate to consider the case
when the basic semiring~$K$ is a quasifield or a semifield.
 This case is
particularly important for the problems of calculus.

\definition{Definition~4.7}
We call a semimodule~$V$ over a quasifield or a semifield~$K$ an {\it
idempotent space}.
\enddefinition

This notion is analogous to the notion of linear (vector) space over a
field.

\definition{Definition~4.8}
We call a semimodule~$V$ over a quasifield~$K$ an {\it idempotent~$a$-space\/} ($b$-{\it space\/}) if it is~$a$-regular (respectively,~$b$-regular) and~$\wh V$ (respectively,~$\wh V_b$) is a
standard semimodule over the~$b$-complete semifield~$\wh K_b$, so
that equality~\thetag{4.7} of Definition~4.4 holds for $x \ne \Bbb I =
\sup V$.
\enddefinition

A quasifield or semifield is an example of a (``one-dimensional'')
idempotent space over itself.

The following statement is a straightforward consequence of
Propositions~3.1 and~3.7, the definitions, and Remark~3.1.

\proclaim{Proposition~4.2}
If~$K$ is a quasifield, then~$\wh K_b$ is an idempotent~$b$-space over~$K$ and~$\wh K_b$, and the semiring~$\wh K$ considered as a
semimodule is an idempotent~$a$-space over~$K$ and~$\wh K_b$ and an~$a$-complete semimodule over~$\wh K$.
\endproclaim

Direct sums and products of semimodules, idempotent spaces,~$a$-spaces, and~$b$-spaces over the same semiring~$K$ can be defined and described in the
usual way and provide a number of new examples.
 Additional nontrivial
examples are generated when subsemimodules and subspaces are considered.

\bhead
4.4.
 Linear maps and functionals
\endbhead
Suppose that~$V$ and~$W$ are idempotent semimodules over an idempotent
semiring~$K$ and {${P\:V\to W}$} is a map from~$V$ to~$W$.  The following
definition is standard.

\definition{Definition~4.9}
 A map $P\:V\to W$
is called {\it additive\/} if
$$
P(x\oplus y)=P(x)\oplus P(y)
\tag4.8
$$
for all $x, y \in V$.
 This map is called {\it homogeneous\/} if
$$
P(k\odot x)=k\odot P(x)
\tag4.9
$$
for all $x\in V$, $k\in K$.
 A map $P\:V\to W$
is called {\it linear\/} if it is additive and homogeneous.
\enddefinition

\definition{Definition~4.10}
Suppose~$V$ and~$W$ are idempotent~$a$-regular ($b$-regular) semimodules
over an idempotent semiring~$K$.
 A linear map $P\:V\to W$ is called~$a$-{\it linear\/}
(respectively,~$b$-{\it linear\/}),
if it is an~$a$-homomorphism (respectively, a~$b$-homomorphism).
\enddefinition

Clearly, the notion of~$a$-linear ($b$-linear) map provides an
algebraic model of the notion of linear (semi)continuous map (see
Sec.~2.8 above; note that in the case of ordinary linear operators
semicontinuity is equivalent to continuity).

\proclaim{Proposition~4.3}
Let a map $P\:V\to W$, where~$V$ and~$W$ are~$a$-regular
\rom($b$-regular\rom) semimodules over a~$b$-regular semiring~$K$,
be~$a$-linear \rom(respectively,~$b$-linear\rom).  Then it is uniquely
extended to a map $\wh P\:\wh V\to \wh W$ \rom(respectively, $\wh P\:\wh
V_b\to \wh W_b$\rom) that is~$a$-linear \rom(respectively,~$b$-linear\rom)
over~$\wh K_b$, where~$\wh V$ and~$\wh W$ are the semimodule completions in
the sense of Definition~4.6.
\endproclaim

\demo{Proof}
The proof is by direct calculation.
  It follows from the definition of~$a$-homomorphism (a~$b$-homomorphism) of idempotent semigroups that the
extension~$\wh P$ is uniquely defined and additive.
 If $k \in \wh K_b$, then $k = \oplus Q$, where $Q \subset K$.
 Obviously,
$$
\wh P(k\odot x)=\wh P(\oplus Q\odot x)
=\oplus P(Q\odot x)=\oplus Q\odot P(x)
=k\odot P(x)=k\odot\wh P(x),
$$
if $x\in V$;
if $x\in \wh V$, then $x=\oplus X$, where
$X\subset V$, and by the
similar argument $\wh P(k \odot x) = k \odot \wh P(x)$.
 This
concludes the proof.
\qed\enddemo

\definition{Definition~4.11}
A {\it functional\/} on a semimodule~$V$ over an idempotent semiring~$K$ is
a map $V \to K$ or a map $V \to \wh K$ from this semimodule to the
completion~$\wh K$ of the semimodule~$K$.
 A functional is called {\it linear\/} if this map is linear.
 A linear functional is called~$a$-{\it
linear\/} ({\it $b$-linear\/}) if it is an~$a$-homomorphism (respectively, a~$b$-homomorphism).
\enddefinition

We assume that an~$a$-linear (a~$b$-linear) functional takes values in
the completion~$\wh K$ (respectively,~$\wh K_b$) if this
completion has the natural structure of semimodule over~$K$.
 This is
always true if~$K$ is a quasifield (by Proposition~3.7 and Remark~3.1).

A general description of~$a$-linear functionals on idempotent spaces is
presented in Sec.~5 below.

\remark{Remark~4.2}
It follows from Proposition~2.4 that a~$b$-linear functional~$f$ on~$V$
is~$a$-linear if and only if $\Up(f(X))=\Up(f(V))$
for any subset~$X$ of~$V$ that is not bounded from above.
 By definition,~$f$ has an extension
$\wh f\:\wh V_b\to \wh K_b$.
 If~$V$ contains the element $\Bbb I=\sup V$, then~$\wh f$ is defined on~$\wh V$ and is~$a$-linear.
 Otherwise~$f$
must be extended to $\wh V=\wh V_b\cup\{\Bbb I\}$,
i.e., defined for the element~$\Bbb I$.
 But $\Bbb I=\sup X$ for any subset~$X$ of~$V$ that is not bounded from above.
 Thus~$\sup f(X)$ must not
depend on the choice of~$X$; then we may put $\wh f(\Bbb I)=\sup
f(X)=\oplus f(X)$.
\endremark

As an example, suppose $\sB(X,\Bbb R)$ is the set of all bounded
functions defined on an arbitrary set~$X$ containing more than one
point; consider $\sB(X,\Bbb R)$ as an idempotent~$b$-space over the
semifield~$\Bbb R(\max, +)$.
The linear functional $\delta_a\:\phi\mapsto\phi(a)$, i.e., the
``delta-function'', is~$b$-linear but not~$a$-linear.  However
a~$b$-linear functional on $\sB(X,\Bbb R)$ always has an~$a$-linear
extension defined on $\sB(X,\wh{\Bbb R})$.

\bhead
4.5.
 Idempotent semimodules and spaces associated with
vector lattices
\endbhead
Examples of idempotent semimodules and spaces that are most important for
Idempotent Analysis are either subsemimodules of (topological) vector
lattices or are dual to these in the sense that they consist of linear
functionals (see Definition~4.11) subject to some regularity conditions,
e.g., of~$a$-linear functionals.

Recall (see~\cite{26, 25, 23}) that a vector space~$V$ over the field of
real numbers~$\Bbb R$ is called {\it ordered\/} (or semiordered~\cite{24})
if~$V$ is equipped with a (partial) order~$\preccurlyeq$ such that all translations
$x \mapsto x + y$ and all homotheties $x \mapsto \lambda x$ preserve this
order (i.e., are isotonic) whenever $\lambda > 0$ and $x, y \in V$.
 If the
operations $x \vee y = x \oplus y = \sup\{x, y\}$ and $x \wedge y =
\inf\{x, y\}$ are defined, then~$V$ is called a {\it vector lattice\/} (in
this case the additive group of~$V$ is lattice ordered, see Sec.~3.3.4
above).
 A {\it topological vector lattice\/} is a Hausdorff topological
vector space~$V$ such that~$V$ is a vector lattice, the lattice operations
$x, y \mapsto x \oplus y$ and $x, y \mapsto x \wedge y$ are continuous, and
the positive cone $\{x\in V\mid x\ge 0\}$ is normal relative
to the topology in~$V$.
In particular, a {\it normed\/} ({\it Banach}) {\it lattice\/} is
a normed (respectively, Banach) space that is a vector lattice such that
the operations $\vee = \oplus$ and $\wedge$ are jointly continuous.

Let~$V$ be a topological vector lattice over~$\Bbb R$, $K$ be a subsemigroup
of the additive group of~$V$ (with respect to the ordinary sum), and~$M$ be a subset of~$V$ invariant under translations $x \mapsto x + k$,
where $k \in K$.
 Suppose further that~$K$ and~$M$ are idempotent
subsemigroups in~$V$ (with respect to $\vee = \oplus$) and the embeddings
$K \to V$ and $M \to V$ are~$b$-homomorphisms.

Obviously,~$V$ is a semifield with respect to the standard idempotent
operations $\oplus = \vee = \sup$ and $\odot = +$ and~$K$ is its
subsemifield.
 Here, the product $K \times M \to M$ is defined by $(k, x)
\mapsto k \odot x = k + x$.

\proclaim{Proposition~4.4}
{\rm1)} The semifields~$V$ and~$K$ are integrally closed \rom(hence, are
quasifields\rom) with respect to the standard idempotent operations $\oplus
= \sup$ and $\odot = +$\rom; in this case $\1 = 0$.

{\rm2)}~$M$ is an idempotent~$b$-space over the quasifield~$K$\rom; the
product $K \times M \to M$ is defined by $k \odot x = k + x$.  In
particular, the whole lattice~$V$ is an idempotent~$b$-space and~$M$ is an
idempotent~$b$-subspace of~$V$ \rom(although it need not be a subspace
of~$V$ in the usual sense\rom).

{\rm3)}
The product $K \times M \to M$ is the restriction of the product
in the~$b$-complete semifield~$\wh V_b$ to $K \times M$.
\endproclaim

\demo{Proof}
We begin with the first statement.
 Let us prove that the semifield~$K$,
which contains an (idempotent) inverse $x^{-1} = -x$ for any of its
elements, is integrally closed.
 Let $x, b \in K$ and $x^n=nx\preccurlyeq b$ for all $n=1,2,\dots$\,.
 Then $x\preccurlyeq(1/n)b$; hence
$x\oplus((1/n)b)=(1/n)b$.
 Using the continuity of the ordinary product and the idempotent
sum, we obtain the following sequence of equalities:
$$
x\oplus\1=x\oplus0
=x\oplus\lim_{n\to\infty}\Bigl(\Bigl(\frac1n\Bigr)b\Bigr)
=\lim_{n\to\infty}\Bigl(x\oplus\Bigl(\frac1n\Bigr)b\Bigr)
=0=\1,
$$
which implies $x\oplus\1= \1$,
i.e., $x\preccurlyeq\1$,
and thus proves the statement.
 The integral closedness of the semifield~$V$ can be shown similarly.
 The rest of Proposition~4.4 can be checked by a direct
calculation.
\qed\enddemo

\definition{Definition~4.12}
We call any idempotent~$b$-space~$M$ described in Proposition~4.4 a {\it
space of $(V, K)$ type}.
\enddefinition

\bhead
4.6.
 Examples
\endbhead

\dhead
\rom{4.6.1}
\enddhead
Suppose~$K$ is a~$b$-complete semifield with zero and~$K$~does not
coincide with the two-element semifield $\{\0, \1\}$.

  The semigroup
$\Map(X, \wh K)$ (see Examples~2.9.4, 3.2.1, 3.6.4) is an idempotent
semimodule over~$K$ (with respect to the pointwise product of functions
by elements of~$K$) but this semimodule is not standard if~$X$ consists of
two or more elements.

The idempotent semimodules $USC(X,\wh{\Bbb R}_{\max})$,
and $LSC(X,\wh{\Bbb R}_{\max})$
of semicontinuous functions are defined similarly; generally, these
semimodules are not standard if~$X$ consists of two or more elements but
they are~$a$-complete.

\dhead
\rom{4.6.2}
\enddhead
The semigroups $USC(X),LSC(X),L^p(X)$, and $\Conv(X, \Bbb R)$ defined in
Examples~2.9.8--2.9.11 (see also Examples~3.2.1 and 3.6) are
idempotent~$a$-spaces over the semifield (and quasifield) $K = \Bbb R(\max,
+)$ under pointwise multiplication of functions by elements of~$K$. The
normal completions of these idempotent semigroups are~$a$-complete
idempotent~$a$-spaces over the~$b$-complete semifield $\Bbb R_{\max}$
and~$a$-complete semimodules over the semiring~$\wh{\Bbb R}_{\max}$.

\dhead
\rom{4.6.3}
\enddhead
The semiring $\wh{\Bbb R}_{\max}$ is an idempotent~$a$-space over the semifield~$\Bbb R_{\max}$.
 Similarly, the set $\wh{\Map(X, \Bbb R)}$ is an idempotent~$a$-space
over~$\Bbb R_{\max}$.
 It is natural to call this space~$n$-dimensional if~$X$
consists of~$n$ elements.

\dhead
\rom{4.6.4}
\enddhead
Any~$a$-complete idempotent semigroup is an idempotent~$a$-space with
respect to the natural action of the Boole algebra, i.e., the semifield
$\{\0, \1\}$.

\dhead
\rom{4.6.5}.
 Spaces of $(V,K)$ type
\enddhead

\dhead
\rom{4.6.5.1}
\enddhead
 The space $C(X)$ of continuous real-valued functions on a
topological space~$X$ (see Examples~2.9.7, 3.2.1, and 3.6.6 above) is a
space of $(V, K)$ type, where~$K$ is the subgroup of constants (i.e.,~$\Bbb R$) and $M = V = C(X)$.

\dhead
\rom{4.6.5.2}
\enddhead
Let~$V = C(X)$, $K$~be the subgroup of integer
constants (i.e.,~$\Bbb Z$, see Example~3.6.2), $M$~be the subgroup of
integer-valued functions.

\dhead
\rom{4.6.5.3}
\enddhead
Let $V = C(\Bbb R)$, $K$~be the subgroup of even
functions, $M = V$.
 Likewise, let $V = C(X)$ and~$K$ be the subgroup of
functions invariant under some group of continuous transformations of~$X$.

\dhead
\rom{4.6.5.4}
\enddhead
Let~$X$ be a measure space with measure~$\mu$.
 Let $V
= M = L^p(X, \mu)$, $K$~be the subgroup of functions invariant under a
group of measure-preserving transformations of~$X$.

\dhead
\rom{4.6.5.5}
\enddhead
Let $$
V=\Bbb R\times\Bbb R=\Bbb R^2, \quad
K=\{(x,x)\mid x\in\Bbb R\}, \quad
M=\bigl\{(x,y)\in V\bigm| |x-y|\le1\bigr\}.
$$
 Note that~$M$ is an ``infinite-dimensional'' subspace of
a two-dimensional space in the sense that it is not a finitely generated
subspace.
 Clearly,~$M$ is not a subspace of~$V$ in the conventional sense.

\ahead
5.
 THE STRUCTURE OF~$a$-LINEAR FUNCTIONALS ON IDEMPOTENT SEMIRINGS
\endahead

\bhead
5.1.
 The basic construction
\endbhead
Suppose~$V$ is an idempotent~$b$-space (e.g., an idempotent~$a$-space,
see Definition~4.8) over a quasifield~$K$.

\definition{Definition~5.1}
Let $x \in \wh V$; by~$x^*$ denote the functional $V \to \wh K$
defined by
$$
y\mapsto x^*(y)
=\wedge\{k\in K\mid y\preccurlyeq k\odot x\}
=\inf_{y\preccurlyeq k\odot x}k,
\tag5.1
$$
for any $y \in V$; we call this functional the~$x$-{\it functional\/}
on~$V$.
\enddefinition

By~$\Bbb I_V$ denote the largest element of~$V$ (if it exists); recall that
$\Bbb I_V=\sup V=\inf\emptyset$,
where~$\emptyset$ is the empty subset
of~$V$.
 Similarly, $\Bbb I_K=\sup K$.
 Note that $\Bbb I^*_V(y)\equiv \0_K$
for all~$y\in V$
(with the possible exception of $y= \Bbb I_V$);
$\0^*_V(y)= \Bbb I_K$, if $y\ne \0_V$,
and $\0^*_V(\0_V)= \0_K$.

\proclaim{Theorem~5.1}
 Let~$V$ be an idempotent~$b$-space over a quasifield~$K$.
 Then for any $x
\in \wh V$ the~$x$-functional $x^*\:y\mapsto x^*(y)$ is an~$a$-linear functional on~$V$.
\endproclaim

\demo{Proof}
Replacing~$K$ and~$V$ by their~$b$-completions and using Definition~4.8
and Proposition~4.1, we may assume that~$V$ is~$b$-complete and~$K$ is a~$b$-complete semifield.
 Then it is clear that $\wh V=V\cup\{\Bbb I_V\}$ and $\wh K=K\cup\{\Bbb I_K\}$.
 For $x= \Bbb I_V$ the statement
of the theorem can be checked directly; consider the case when $x\prec\Bbb I_V$.
 The fact that~$V$ is standard (see Definitions~4.4 and 4.8) implies
that $y\preccurlyeq x^*(y)\odot x$,
if the set $K(y)=\{k\in K\mid y\preccurlyeq k\odot x\}$
is not empty.

Let~$Y$ be an arbitrary subset of the semimodule~$V$.
 Using the
construction of the map $y \mapsto x^*(y)$, it is easy to prove that this
map preserves the order (i.e., is isotonic); therefore $\oplus x^*(Y)
=\sup\{x^*(y)\mid y\in Y\}\preccurlyeq x^*(\oplus Y)$.

If the set $K(y)=\{k\in K\mid y\preccurlyeq k\odot x\}$ is empty, then
$$
x^*(\oplus Y)= \Bbb I_K,\quad
\oplus x^*(Y)=\sup\{x^*(y)\mid y\in Y\}= \Bbb I_K,
$$
where $\Bbb I_K=\sup K=\sup\wh K$; thus in this case $x^*(\oplus Y)=\oplus
x^*(Y)$. Otherwise,
$$
(\oplus x^*(Y))\odot x\succcurlyeq x^*(y)\odot x\succcurlyeq y
\text{ for any } y\in Y.
$$
Thus $(\oplus x^*(Y))\odot
 x\succcurlyeq\oplus Y$; therefore, $\oplus x^*(Y)\succcurlyeq x^*(\oplus
Y)$.  Since we have already proved that $\oplus x^*(Y)\preccurlyeq
x^*(\oplus Y)$, it follows that
$$
x^*(\oplus Y)=\oplus x^*(Y),
\tag5.2
$$
i.e., the functional $y \mapsto x^*(y)$ is an~$a$-homomorphism.

Now let us prove that this functional is homogeneous, i.e., that
$$
k\odot x^*(y)=x^*(k\odot y)
\tag5.3
$$
for all $k\in K$ and $y\in V$.
 Suppose that~$p$ is an invertible element
of~$K$ and~$y$ is any element of~$V$.  Since multiplication by~$p$, or homothety, is an automorphism of~$K$, we see that $$
p\odot K(y)
=p\odot\{k\in K\mid y\preccurlyeq k\odot x\}
=\{k\in K\mid p\odot y\preccurlyeq k\odot x\}
=K(p\odot y).
$$
 Therefore the set $K(p \odot y)$ is empty
if the set $K(y)$ is empty; thus $x^*(y)= x^*(p\odot y)= \Bbb I_K$ and
$p\odot x^*(y)=x^*(p\odot y)$ since $k\odot\Bbb I_K= \Bbb I_K$.
 Now if the set
$K(y)$ is not empty, then $K(p \odot y)$ is not empty and
$$
p\odot x^*(y)
=p\odot\wedge K(y)=\wedge(p\odot K(y))
=\wedge K(p\odot y)=x^*(p\odot y).
$$
 Hence the homogeneity~\thetag{5.3} is proved for all invertible, i.e.,
nonzero, elements $p = k \in K$.
 Finally, if $k = \0$, then $\0\odot x^*(y)=\0=x^*(\0)=x^*(\0\odot y)$.
 We have proved that the functional is
homogeneous; this completes the proof of Theorem~5.1.
\qed\enddemo

\bhead
5.2.
 The basic theorem on the structure of functionals
\endbhead

\proclaim{Proposition~5.1}
Suppose~$V$ is a standard semimodule over a~$b$-complete semiring~$K$
\rom(e.g.,~$V$ is a~$b$-complete~$b$-space, see Definitions~4.4 and
4.8\rom{);} let~$f$ be an~$a$-linear functional on~$V$ such that~$f$ takes
the value~$\1$ and $f(\Bbb I_V)\succ\1$, where $\Bbb I_V=\sup V$. Then
there exists a unique~$x$-functional $y \mapsto x^*(y)$ such that
$x^*(y)\succcurlyeq f(y)$ for any $y\in V$ and $x^*(y)=f(y)$ if $f(y)$ is
invertible in~$K$\rom; here $x=\oplus\{y\in V\mid f(y)\preccurlyeq\1\}$.
\endproclaim

\demo{Proof}
Let $x=\oplus\{y\in V\mid f(y)\preccurlyeq\1\}$; clearly, $f(x)= \1$.
Since by the hypothesis $f(\Bbb I_V)\succ\1$, we see that $x\ne \Bbb I_V$.
If $y\in V$ and $f(y)$ is invertible, then $\1=(f(y))^{-1}\odot
f(y)=f((f(y))^{-1}\odot y)$; thus $x\succcurlyeq(f(y))^{-1}\odot y$, so
$f(y)\odot x\succcurlyeq y$ and therefore $x^*(y)\preccurlyeq f(y)$. On the
other hand, since $x^*(y)\odot x\succcurlyeq y$, it follows that
$x^*(y)=x^*(y)\odot f(x)\succcurlyeq f(y)$.  As a result, $x^*(y)=f(y)$,
if~$f(y)$ is invertible.
\qed\enddemo

\proclaim{Lemma~5.1}
Let~$K$ be a~$b$-complete semifield that does not coincide with the Boolean
algebra $\{\0, \1\}$.
 Then
$$
\0=\wedge(K\setminus\{\0\})=\inf(K\setminus\{\0\}).
$$
\endproclaim

\demo{Proof}
Let $m = \wedge(K \setminus \{\0\})$.
 It follows from the assumption that
there exists an invertible element $k \ne \1$; we may assume that $k \prec \1$ (in the converse case replace $k$ by $(\1 \oplus k)^{-1}$).
 Thus it follows from the definition of~$m$ that $m\preccurlyeq k\prec\1$.
 If $m \ne \0$, then the elements~$m$ and $m \odot m$ are invertible and $m \odot m
\prec m \odot \1 \prec m$, which contradicts the definition of the element~$m$. This concludes the proof of the lemma.
\qed\enddemo

\proclaim{Theorem~5.2}
Suppose~$V$ is an idempotent~$b$-space over a quasifield~$K$\rom; then any
nonzero~$a$-linear functional~$f$ defined on~$V$ has the form $f = x^*$ for
a unique $x \in \wh V_b$.  If $K\ne\{\0,\1\}$, then
$$
x=\oplus\{y\in V\mid f(y)\preccurlyeq\1\}.
$$
In the converse case $x=\oplus\{y\in V\mid f(y)= \0\}$.
\endproclaim

\demo{Proof}
Replacing~$K$ and~$V$ by their~$b$-completions and using Definition~4.8
and Proposition~4.1, we may assume that~$V$ is~$b$-complete and~$K$ is a~$b$-complete semifield.
 Suppose~$f$ is an arbitrary nonzero~$a$-linear
functional defined on~$V$ and $\wh f \: \wh V \to \wh K$
is its extension to~$\wh V$ according to the definitions of~$a$-linearity and~$a$-homomorphism (see Definitions~4.10 and 2.5-a).

First suppose that $K = \{\0, \1\}$; let $x=\oplus\{y\in V\mid f(y)= \0\}$.
It is easily shown that in this case $f(y) = \0$ if $y\preccurlyeq x$ and
$f(y)= \1$ otherwise; thus $f = x^*$.  If $x= {\Bbb I}_V$, then $\wh
f({\Bbb I}_V)= \0$, which is impossible if~$f$ is nonzero; thus $x\in V=
\wh V_b$.

Now let $K\ne\{\0,\1\}$.
Lemma~5.1 implies that $\0 = \wedge
\{K \setminus \{\0\})$.
 Let $x=\oplus\{y\in V\mid f(y)\preccurlyeq\1\}$.
 If $x= \Bbb I_V$, then $\wh f(\Bbb I_V)\preccurlyeq\1$,
which is impossible
if~$f$ is nonzero since by homogeneity the range of~$f$ must contain all
nonzero (i.e., invertible) elements of~$K$, including some $k \succ \1$.
Thus $x \in V = \wh V_b$ and, evidently, $f(x) = \1$.

In the proof of Proposition~5.1 it was shown that $x^*\succcurlyeq f$ and $x^*(y)=f(y)$ if $f(y)\in K\setminus\{\0\}$.
 Thus in order to prove that $x^* = f$ it is sufficient to check that $x^*(y) = \0$ if $f(y) = \0$.
 Indeed, if
$f(y) = \0$, then $f(k \odot y) = \0 \prec \1$ if $k\ne \0$.
 But this implies that $k\odot y\preccurlyeq x$ if $k\ne \0$; therefore
$k\odot x\succcurlyeq y$ if $k\ne \0$, so $x^*(y)\preccurlyeq\wedge(K\setminus\{\0\})
=\0\preccurlyeq x^*(y)$, which completes the proof.
\qed\enddemo

\remark{Remark~5.1}
Note that earlier theorems of this type (see, e.g.,~\cite{7, 8, 11,~14})
were proved using the assumption that the semimodule contains
sufficiently many ``Dirac delta functions''.
 On the contrary, our proof of
Theorem~5.2 is based on the construction of ``lower envelopes''.
 Thus
Theorem~5.2 is valid in the cases when ``delta functions'' do not exist (e.g.,
in spaces of integrable functions).
 Hence the area of application of
Theorem~5.2 is to a large degree complementary to that of the results based
on ``delta functions''.
  Moreover, in general the property of having a
sufficient collection of ``delta functions'' is not inherited by
subsemimodules.
\endremark

\bhead
5.3.
 Theorems of the Hahn--Banach type
\endbhead
Let~$V$ be a semimodule over an idempotent semiring~$K$.
 A {\it subsemimodule}~$W$ of the semimodule~$V$ is a subsemigroup~$W$ in~$V$ that is closed under multiplication by elements (coefficients) of~$K$.
 Note that~$W$ itself is a semimodule over~$K$.

\definition{Definition~5.2}
Let~$V$ be an idempotent~$a$-space ($b$-space) over a quasifield~$K$.
 A
subsemimodule~$W$ of~$V$ is called an~$a$-{\it subspace\/} (respectively, a~$b$-{\it subspace\/}) of~$V$ if the embedding $i \: W \to V$ has a unique~$a$-linear (respectively,~$b$-linear) extension $\wh W \to \wh V$ (respectively, $\wh W_b \to \wh V_b$) to the
completions of the semimodules defined over~$\wh K_b$.
\enddefinition

The basic definitions immediately imply the following statement.

\proclaim{Proposition~5.2}
Let~$V$ be an idempotent~$a$-space \rom($b$-space\rom) over a quasifield~$K$ and~$W$ be its~$a$-subspace \rom(respectively,~$b$-subspace\rom).
 Then~$W$ and its completion~$\wh W$ \rom(respectively,~$\wh W_b$\rom) are idempotent~$a$-spaces \rom(respectively,~$b$-spaces\rom).
\endproclaim

\proclaim{Theorem~5.3}
Let~$V$ be an idempotent~$b$-space over a quasifield~$K$ and~$W$ be its~$b$-subspace. Then any~$a$-linear functional on~$W$ has an~$a$-linear extension to~$V$.
\endproclaim

This result is an immediate consequence of Theorem~5.2 and Proposition~5.2.

\proclaim{Theorem~5.4}
Let~$V$ be an idempotent~$b$-space.
 If $x, y \in V$ and $x\ne y$, then there exists an~$a$-linear functional~$f$ on~$V$ such that $f(x)\ne f(y)$.
\endproclaim

\demo{Proof}
If $x \succ y$, then $y^*(x) \succ \1$ and $y^*(y)\preccurlyeq\1$, so $f= y^*$ has the required property.
 Indeed, in the converse case we see that $x\prec\Bbb I_V$,
the inequality $x^*(y)\preccurlyeq\1$
is not possible, but $x^*(x)\preccurlyeq\1$; thus $f= x^*$ has the required property.
\qed\enddemo

\remark{Remark~5.2}
It follows from the definition of spaces of $(V, K)$ type (see Sec.~4.6 above) that their $b$-subspaces are spaces of
$(V, K)$ type themselves.
\endremark

\remark{Remark~5.3}
The subsemimodule $C(X)$ in $USC(X)$ is not a~$b$-subspace since the
corresponding embedding is not a~$b$-homomorphism. Still, $C(X)$ is an
idempotent~$b$-space (and~$a$-space) over $\Bbb R(\max, +)$, and~$a$-linear
functionals on $C(X)$ have~$a$-linear extensions to $USC(X)$. It is easy to
extend Theorem~5.3 to cases of this kind.
\endremark

\bhead
5.4.
 Analogs of the Banach--Steinhaus theorem and the closed graph
theorem
\endbhead
The following statements are straightforward consequences of the definitions;
they can be considered as analogs of well-known results of the
traditional Functional Analysis (the Banach--Steinhaus theorem and the
closed graph theorem).  {\it In this section all~$a$-spaces are assumed to be~$a$-complete.} The results can easily be extended to the case of
incomplete spaces by means of the completion procedure.\footnote{An
idempotent analog of the Banach--Steinhaus theorem for spaces of continuous
functions that tend to $\0 = \infty$ at infinity is formulated in~\cite{7, p.~59}; see also~\cite{8, p.~52}.
 These spaces are not~$a$-complete but can be completed in the usual way.}

\proclaim{Proposition~5.3}
Suppose $\sP=\{P_\alpha\}$ is a family of~$a$-linear maps of
an $a$-space~$V$ to an $a$-space~$W$ and the map $P \: V \to W$ is the
pointwise sum of these maps, i.e., $P(x)=\sup\{P_\alpha(x)\mid P_\alpha\in
\sP\}$\rom; then the map~$P$ is~$a$-linear.
\endproclaim

Indeed, we have
$$
\align
P(k\odot x)
&=\sup\{P_\alpha(k\odot x)\mid P_\alpha\in\sP\}
=\sup\{k\odot P_\alpha(x)\mid P_\alpha\in\sP\}
\\
&=k\odot\sup\{P_\alpha(x)\mid P_\alpha\in\sP\}
=k\odot P(x),
\endalign
$$
for any $x \in V$; thus the map~$P$ is homogeneous. If $X \subset V$,
then $$
\align
P(\oplus X)
&=\sup\{P_\alpha(\oplus X)\mid P_\alpha\in\sP\}
=\sup\{\oplus P_\alpha(X)\mid P_\alpha\in\sP\}
\\
&=\sup\{P_\alpha(x)\mid x\in X,\ P_\alpha\in\sP\}
=\sup\bigl\{\sup\{P_\alpha(x)\mid x\in X\}\mid
P_\alpha\in\sP\bigr\}
\\
&=\sup\{P(x)\mid x\in X\}
=\oplus P(X);
\endalign
$$
this completes the proof.\qed

\proclaim{Corollary}
The pointwise sum \rom(i.e., the least upper bound\rom) of a family of~$a$-linear functionals is an~$a$-linear functional.
\endproclaim

\proclaim{Proposition~5.4}
Let~$V$ and~$W$ be~$a$-spaces.
 A linear map $P \: V \to W$ is~$a$-linear if and only if its graph~$\Gamma$ in $V \times W$ is closed under
sums \rom(i.e., least upper bounds\rom) of arbitrary subsets.
\endproclaim

The hypothesis of Proposition~5.4 implies that the embedding $i \: \Gamma \to V \times W$ is~$a$-linear.
 To complete the proof, it is
sufficient to note that~$P$ is the composition of three~$a$-linear
maps: the isomorphism $x \mapsto (x, P(x)) \in \Gamma$, the embedding~$i$,
and the  projection $V \times W \to W$.\qed

\bhead
5.5.
 The scalar product
\endbhead
Let us consider a class of idempotent spaces on which a natural structure of
scalar product can be defined.

\definition{Definition~5.3}
An idempotent~$b$-space~$A$ over a quasifield~$K$ is called an {\it
idempotent $b$-semi\-alge\-bra\/} over~$K$ if~$A$ is equipped with the
structure of idempotent semiring consistent with the product $K \times A
\to A$ (in the sense that the product is associative). If the product $A
\times A \to A$ is a separate~$a$-homomorphism ($b$-homomorphism), then
the~$b$-semialgebra~$A$ is called~$a$-{\it regular} (respectively,~$b$-{\it
regular}).
\enddefinition

The theory of vector lattices is an important source of examples of
idempotent~$b$-semialgebras.

\proclaim{Proposition~5.5}
For any invertible element $x \in A$, where~$A$ is a~$b$-semialgebra, and
any element $y \in A$, the following equality holds:
$$
x^*(y)=\1^*_A(y\odot x^{-1}).
\tag5.4
$$
\endproclaim

\demo{Proof}
The completion procedure allows to reduce the proof to the case in
which~$K$ is a~$b$-complete semifield.  In this case
$$
x^*(y)=\inf\{k\mid k\odot x\succcurlyeq y\}
=\inf\{k\mid(k\odot x)\odot x^{-1}\succcurlyeq y\odot x^{-1}\}
=\inf\{k\mid k\odot\1_A\succcurlyeq y\odot x^{-1}\}.
$$
 To conclude the proof, note that $\inf\{k\mid k\odot\1_A\succcurlyeq y\odot x^{-1}\}
=\1^*_A(y\odot x^{-1})$.
\qed\enddemo

\definition{Definition~5.4}
Suppose~$A$ is a commutative~$b$-semialgebra over a quasifield~$K$.
 The
map $A \times A \to \wh K$ defined by the formula $(x,y)\mapsto\<x,y\>=\1^*(x\odot y)$ is called the {\it canonical scalar
product\/} (or simply the {\it scalar product\/}) on~$A$.
\enddefinition

\proclaim{Proposition~5.6}
The scalar product $\<x,y\>$ on a commutative~$a$-regular
\rom($b$-regular\rom)~$b$-semialgebra over a quasifield~$K$ has the
following properties\rom:
\roster
\item"\rom{1)}" the map $(x, y) \mapsto \< x, y\>$ is a
separate~$a$-homomorphism \rom(respectively,~$b$-homomorphism\rom{);}
\item"\rom{2)}" for all $k\in K$, $x,y\in A$ and any \rom(respectively, any
bounded from above\rom) subset~$X$ of~$A$, the following equalities hold\rom:
$$
\align
\<x,y\>
&=\<y,x\>,
\tag5.5
\\
\<k\odot x,y\>
&=\<x,k\odot y\>=k\odot\<x,y\>,
\tag5.6
\\
\<\oplus X,y\>
&=\bigoplus_{x\in X}\<x,y\>.
\tag5.7
\endalign
$$
\endroster
\endproclaim

Proposition~5.6 follows directly from the definitions; formula~\thetag{5.7}
follows from the first statement.

The following statement is an idempotent analog of the well-known
Riesz--Fischer theorem.

\proclaim{Theorem~5.5}
Let~$A$ be a commutative~$b$-semialgebra over a quasifield~$K$ such that~$A$ itself is a quasifield.
 Then any nonzero~$a$-linear functional~$f$
on~$A$ has the form
$$
f(y)=\<x,y\>,
\tag5.8
$$
where $x\in \wh A_b$, $x\ne \0$, and
$\<\,\cdot\,,\,\cdot\,\>$
is the canonical scalar product on~$\wh A_b$.
\endproclaim

This theorem follows directly from Remark~3.1, Theorem~5.2, and
Proposition~5.5.

\remark{Remark~5.4}
Evidently, the canonical scalar product can be defined for the case in
which the~$b$-semialgebra~$A$ is not commutative, and Theorem~5.5 can be
extended to this case as well.  \endremark


\bhead
5.6.
 The skew-scalar product and elements of duality
\endbhead

\dhead
\rom{5.6.1}
\enddhead
Let~$V$ be an idempotent~$b$-space over a~$b$-complete
semifield~$K$.
 For any two elements $x, y \in V$, by $[x, y]$ denote the
value $x^*(y)$ of the functional~$x^*$ at the element~$y$ (see
Definition~5.1).

\definition{Definition~5.5}
We say that the  map $V \times V \to \wh K$ defined by $(x, y) \mapsto [x, y] = x^*(y)$ is the {\it canonical skew-scalar product\/} (or simply the
{\it skew-scalar product\/}) on~$V$.
\enddefinition

The following statements are direct consequences of our definitions.

\proclaim{Proposition~5.7}
The skew-scalar product has the following properties:
$$
\gather
[x,x]\preccurlyeq\1,
\tag5.9
\\
[x,k_1\odot y_1\oplus k_2\odot y_2]
=k_1\odot[x,y_1]\oplus k_2\odot[x,y_2],
\tag5.10
\\
[k\odot x,y]=k^{-1}\odot[x,y],
\tag5.11
\\
[x_1\wedge x_2,y]=[x_1,y]\oplus[x_2,y]
\tag5.12
\endgather
$$
for all $x,x_1,x_2\in \wh V$, $y,y_1,y_2\in V$,
$k,k_1,k_2\in K$, $k\ne \0$.
\endproclaim

\proclaim{Proposition~5.8}
If the canonical scalar product $(x,y)\mapsto\<x,y\>$ is
defined on a space~$V$ \rom(i.e.,~$V$ is a~$b$-semialgebra, see Sec.~5.5
above\rom), then the following equalities hold\rom:
$$
\gather
\<x,y\>=[y^{-1},x],
\qquad
[x,y]=\<x^{-1},y\>,
\tag5.13
\\
[x,y]=[y^{-1},x^{-1}]
\tag5.14
\endgather
$$
for all invertible $x,y\in V$.
\endproclaim

\dhead
\rom{5.6.2}
\enddhead
Suppose~$V$ and~$W$ are~$a$-complete~$a$-spaces over a~$b$-complete
semifield~$K$ that does not coincide with the Boolean algebra $\{\0, \1\}$
and~$W$ is an~$a$-subspace of~$V$.

By~$V^*$ denote the set of all~$a$-linear functionals on~$V$;
this set
forms an idempotent semimodule under pointwise operations.
 It
follows from Theorem~5.2 that the sets of elements of~$V$ and~$V^*$ are in
one-to-one correspondence; however the respective structures of semimodules
over~$K$ are different.

The following result follows from Definitions~4.3, 4.4, 4.8, and 4.11,
Propositions~3.1 and 5.7, and Theorems~5.1 and 5.2.

\proclaim{Theorem~5.6}
A functional $y \mapsto f(y)$ on an~$a$-complete~$a$-space~$V$ is~$a$-linear if and only if it has the form $f(y) = [x, y] = x^*(y)$, where
$x \in V$.
 An idempotent semimodule~$V^*$ over~$K$ is an~$a$-complete~$a$-space.
 In addition,
$$
x^*_1\oplus x^*_2=(x_1\wedge x_2)^*,
\quad
k\odot x^*=(k^{-1}\odot x)^*,
\quad
\0^*_V=\Bbb I_{V^*}=\sup V^*,
\quad
\Bbb I^*_V=\sup V=\0_{V^*}
$$
and the canonical order on~$V^*$ is opposite to that on~$V$.
\endproclaim

Define a map~$P \: V \to W$ of the space~$V$ to its~$a$-subspace~$W$ by the formula
$$
P(x)=\inf\{w\in W\mid w\succcurlyeq x\}.
\tag5.15
$$

\proclaim{Proposition~5.9}
The map~$P$ is an~$a$-linear projection.
\endproclaim

The proof of this statement is similar to that of Theorem~5.1.

\proclaim{Proposition~5.10}
The subspace~$W$ is the set of all solutions to the system of equations
$$
[y,x]=[y,P(x)],
$$
where~$y$ runs over~$V$ and the projection~$P$ is defined
by~\thetag{5.15}.
\endproclaim

This statement is a straightforward consequence of Theorem~5.4.

\proclaim{Theorem~5.7}
The map $x \mapsto x^{**} = (x^*)^*$ is an isomorphism of~$a$-spaces $V \to V^{**} = (V^*)^*$.
\endproclaim

This theorem follows from Theorem~5.6 since it can easily be checked that
$x^{**}(y^*) = y^*(x)$.

\bhead
5.7.
 Examples
\endbhead
The semifield $\sB(X)$
of all bounded real-valued functions on an
arbitrary set~$X$ (see Examples~2.9.6 and 3.2.1) is a~$b$-semialgebra
over the semiring $K = \Bbb R(\max, +)$.

In this case,
$$
\1^*(\phi)
=\sup_{x\in X}\phi(x)
=\int^\oplus_X\phi(x)\,dx
\tag5.16
$$
and the scalar product can be expressed in terms of idempotent integration
(see the introduction, Sec.~1.2):
$$
\<\phi_1,\phi_2\>
=\sup_{x\in X}\bigl(\phi_1(x)\odot\phi_2(x)\bigr)
=\int^\oplus_X\bigl(\phi_1(x)+\phi_2(x)\bigr)\,dx,
\tag5.17
$$
where $\phi_1,\phi_2\in\sB(X)$.

Scalar products similar of the~\thetag{5.17} type have been systematically
used in Idempotent Analysis (see, e.g.,~\cite{1--11, 14, 42--43,~49})
in specific spaces, while $a$-linear functionals on idempotent spaces
(including spaces of $(V, K)$ type defined in Sec.~4.5) can easily be
described in terms of idempotent measures and integrals by means of
Theorems~5.2, 5.5, and~5.6.

For example, the idempotent semiring of convex functions $\Conv(X, \Bbb R)$
is a $b$-subspace of the $b$-semialgebra $\Map(X, \Bbb R)$ over the
semifield (and quasifield) $K = \Bbb R(\max, +)$.

Any nonzero~$a$-linear functional~$f$ on $\Conv(X)$ has the form
$$
\phi\mapsto f(\phi)
=\sup_x\{\phi(x)+\psi(x)\}
=\int^\oplus_X\phi(x)\odot\psi(x)\,dx,
$$
where~$\psi$ is a concave function, i.e., an element of the idempotent
semiring
$$
\Conc(X, \Bbb R) = -\Conv(X, \Bbb R) \subset \Map(X, \Bbb R),
$$
see Examples~2.9.9, 2.9.10, and~3.2.1.

\ahead
6.
 COMMENTARY
\endahead

\chead
6.1
\endchead
When this paper was already finished, V.~N.~Kolokoltsov kindly called the
authors' attention to several early papers~\cite{55--59}
and~A.~N.~Sobolevski\u{\i} pointed to the paper~\cite{36}.
In~\cite{55--58,~36}, elements of matrix and linear algebra over
idempotent semirings were considered, including the eigenvector and
eigenvalue problem~\cite{57}.  In~\cite{36}, the matrix calculus over $\Bbb
R_{\min}$ was applied to optimization problems on graphs.  In the
remarkable introduction to his paper~\cite{57}, {N.~N.~Vorob{\mz}ev}
predicted the onset of Idempotent Mathematics as a far-reaching new field
of mathematics.  N.~N.~Vorob{\mz}ev\- called Idempotent Mathematics
`extremal mathematics' and idempotent semimodules `extremal spaces'.
Developing N.~N.~Vorob{\mz}ev's\- ideas, A.~A.~Korbut announced
in~\cite{59} several results including a variant of the Hahn--Banach
theorem (on linearity under finite combinations of elements) and a
finite-dimensional analog of the Riesz--Fischer theorem.  In the same
framework, K.~Zimmermann~\cite{60} presented a very general geometrical
variant of the Hahn--Banach theorem.  Recently, G.~Cohen, S.~Gaubert, and
J.-P.~Quadrat proved a separation theorem of the Hahn--Banach type for a
point and a subsemimodule in the case of the finite-dimensional free
semimodule $\wh{\Bbb R}^n_{\max}$ over the semiring $\wh{\Bbb R}_{\max}$
and described analogs of the basic notions of ``Euclidean'' geometry (see
the paper~\cite{61}, which was kindly sent to us by its authors).  Note that
the earliest paper on idempotent linear algebra that we know (with
applications to the theory of finite automata) belongs to
S.~Kleene~\cite{15}.  For a good survey of idempotent linear algebra,
see~\cite{62}; see also~\cite{6--8, 11--13, 16--21}.

The first stage of development of Idempotent Analysis is presented in the
books~\cite{1, 2, 4} and the papers~\cite{3, 5}.
 The next stage, at which
important results of Idempotent Analysis were obtained for some specific
functional spaces, is represented by the books~\cite{6--8}.
V.~N.~Kolokoltsov's survey~\cite{63} describes aspects of applications of
Idempotent Analysis.

In many papers (see, e.g.,~\cite{64--66}), elements of Idempotent
Mathematics appeared implicitly.
 For additional historical information and
references see~\cite{8, 9, 12, 13, 19--21, 49, 52,~53}.

\chead
6.2
\endchead
Idempotent quantization (dequantization), which was discussed in the
introduction, is related to logarithmic transformations that date back to
the classical papers of E.~Schr\"{o}dinger~\cite{67} and E.~Hopf~\cite{68}.
The subsequent progress of E.~Hopf's ideas has culminated in the vanishing
viscosity method (the method of viscosity solutions).
 This method was
developed by P.-L.~Lions and others for the study of solutions to first and
second order PDE and the related problems~\cite{69--73}.
 A different
approach, which leads to similar results, was considered in the framework
of Idempotent Analysis (see, e.g.,~\cite{3, 6--8}).

The Lax--Ole\u{\i}nik formula was introduced in~\cite{74, 75}; see also an
important survey paper by A.~I.~Subbotin~\cite{72}.

Recently O.~Viro~\cite{76} described the result of idempotent
dequantization of real algebraic geometry.

\chead
6.3
\endchead
The main results presented in Sec.~5 of this paper were announced
in~\cite{77}, which was preceded by preliminary publications~\cite{78,~79}.
The algebraic approach was applied in~\cite{80} to the construction of an
idempotent version of topological tensor products and kernel operators in
the spirit of Grothendieck~\cite{81}.
 Papers devoted to idempotent kernel
spaces, eigenvector theorems, an idempotent version of the Schauder
theorem, and the idempotent representation theory, with applications to the
idempotent harmonic analysis, are now in preparation.


\ahead 
ACKNOWLEDGMENTS
\endahead


The authors wish to express their thanks to~V.~N.~Kolokoltsov for useful
comments; they are especially grateful to~A.~N.~Sobolevski\u{\i} who
read the text carefully and made a number of corrections.

The English translation was thoroughly examined and corrected by
A.~B.~Sossinsky.

This research was supported by INTAS and the Russian Foundation for Basic
Research (RFBR) under the joint INTAS--RFBR grant no.~95-91 and the RFBR
grant no.~99-01-00196.

\par\removelastskip\penalty-200\medskip

\Refs


\item{1.}
S.~M.~Avdoshin, V.~V.~Belov, and V.~P.~Maslov,
{\it Mathematical Aspects of Computational Media Synthesis\/}
[in Russian],
MIEM,
Moscow,
1984.

\item{2.}
V.~P.~Maslov,
{\it Asymptotic Methods for Solving Pseudodifferential Equations\/}
[in Russian],
Nauka,
Moscow,
1987.

\item{3.}
V.~P.~Maslov,
``On a new superposition principle for optmization problems,''
{\it Uspekhi Mat. Nauk\/} [{\it Russian Math. Surveys\/}],
{\bf 42}
(1987),
no.~3,
39--48.

\item{4.}
V.~P.~Maslov,
{\it M\'ethodes op\'eratorielles},
Mir,
Moscow,
1987.

\item{5.}
S.~M.~Avdoshin, V.~V.~Belov, V.~P.~Maslov, and A.~M.~Chebotarev,
``Design of computational media: mathematical aspects,''
in: {\it Mathematical Aspects of Computer Engineering}
(V.~P.~Maslov, K.~A.~Volo\-sov, editors),
Mir,
Moscow,
1988,
pp.~9--145.

\item{6.}
{\sl Idempotent Analysis}
(V.~P.~Maslov, S.~N.~Samborski\u\i, editors),
vol.~13,
Adv. Sov. Math,
Amer. Math. Soc.,
Providence (R.I.),
1992.

\item{7.}
V.~P.~Maslov and V.~N.~Kolokol{\mz}tsov,
{\it Idempotent Analysis and Its Application to Optimal Control\/}
[in Russian],
Nauka,
Moscow,
1994.

\item{8.}
V.~N.~Kolokoltsov and V.~P.~Maslov,
{\it Idempotent Analysis and Applications},
Kluwer Acad. Publ.,
Dordrecht,
1997.

\item{9.}
G.~L.~Litvinov and V.~P.~Maslov,
{\it Correspondence Principle
for Idempotent Calculus and Some Computer
Applications}
(IHES/M/95/33),
Institut des Hautes \'Etudes Scientifiques,
Bures-sur-Yvette,
1995,
in: {\it Idempotency},
Publ. of the Newton Institute,
Cambridge Univ. Press,
Cambridge,
1998,
pp.~420--443.

\item{10.}
G.~L.~Litvinov and V.~P.~Maslov,
``Idempotent mathematics: the correspondence principle and its applications
to computing,''
{\it Uspekhi Mat. Nauk\/} [{\it Russian Math. Surveys\/}],
{\bf 51}
(1996),
no.~6,
209--210.

\item{11.}
P.~I.~Dudnikov and S.~N.~Samborskii,
``Endomorphisms of semimodules over semirings with idempotent operation,''
{\it Izv. Akad. Nauk SSSR Ser. Mat.} [{\it Math. USSR-Izv.}],
{\bf 55}
(1991),
no.~1,
91--105.

\item{12.}
F.~L.~Bacelli, G.~Cohen, G.~J.~Olsder, and J.-P.~Quadrat,
{\it Synchronization and Linearity: an Algebra
for Discrete Event Systems},
John Wiley \& Sons Publ.,
New York,
1992.

\item{13.}
{\sl Idempotency}
(J.~Gunawardena, editor),
Publ. of the Newton Institute,
Cambridge Univ. Press,
Cambridge,
1998.

\item{14.}
M.~A.~Shubin,
{\it Algebraic Remarks on Idempotent Semirings and the Kernel Theorem
in Spaces of Bounded Functions\/}
[in Russian],
Institute for New Technologies,
Moscow,
1990.

\item{15.}
S.~C.~Kleene,
``Representation of events
in nerve nets and finite automata,''
in: {\it Automata Studies}
(J.~McCarthy, C.~Shannon, editors),
Princeton Univ. Press,
Princeton,
1956,
pp.~3--40.

\item{16.}
B.~A.~Carr\'e,
``An algebra
for network routing problems,''
{\it J. Inst. Math. Appl.},
{\bf 7}
(1971),
273--294.

\item{17.}
B.~A.~Carr\'e,
{\it Graphs and Networks},
The Clarendon Press, Oxford Univ. Press,
Oxford,
1979.

\item{18.}
M.~Gondran and M.~Minoux,
{\it Graphes et algorithmes},
Eyrolles,
Paris,
1979;
English translation of the newer edition:
{\it Graphs and Algorithms},
John Wiley and Sons Publ.,
New York
1984.

\item{19.}
R.~A.~Cuningham-Green,
{\it Minimax Algebra},
vol.~166,
Springer Lect. Notes
in Economics and Math. Systems,
1979.

\item{20.}
U.~Zimmermann,
{\it Linear and Combinatorial Optimization
in Ordered Algebraic
Structures},
vol.~10,
Ann. Discrete Math,
1981.

\item{21.}
J.~S.~Golan,
{\it Semirings and Their Applications},
Kluwer Acad. Publ.,
Dordrecht,
1999.

\item{22.}
{\sl General Algebra}
(L.~A.~Skornyakov, editor)
[in Russian],
Mathematics Reference Library,
vol.~1,
Nauka,
Moscow,
1990;
vol.~2,
ibid.,
1991.

\item{23.}
G.~Birkhoff,
{\it Lattice Theory},
Amer. Math. Soc.,
Providence,
1967.

\item{24.}
B.~Z.~Vulikh,
{\it Introduction to the theory of semiordered spaces\/}
[in Russian],
Fizmatgiz,
Moscow,
1961.

\item{25.}
M.~M.~Day,
{\it Normed Linear Spaces},
Springer,
Berlin,
1958.

\item{26.}
H.~H.~Schaefer,
{\it Topological vector spaces},
the Macmillan Company,
New York, and
Collier--Macmillan Ltd.,
London,
1966.

\item{27.}
L.~Fuchs,
{\it Partially Ordered Algebraic Systems},
Pergamon Press,
Oxford etc.,
1963.

\item{28.}
P.~Del~Moral,
``A survey of Maslov optimization theory,''
in: {\it Idempotent Analysis and Applications},
Kluwer Acad. Publ.,
Dordrecht,
1997.

\item{29.}
T.~Huillet, G.~Rigal, and G.~Salut,
{\it Optimal Versus Random Processes: a General Framework},
IFAC World Congress, Tallin, USSR, 17--18~August~1990,
CNRS-LAAS Report no.~242~89251, juillet~1989, 6p.,
Toulouse,
1989.

\item{30.}
{\sl Proceedings of the 11th Conf. on Analysis and Optimization of
Systems: Discrete Event Systems}
(G.~Gohen, J.-P.~Quadrat, editors),
vol.~199,
Lect. Notes in Control and Information Sciences,
Springer,
1994.

\item{31.}
P.~Del~Moral, J.-Ch.~Noyer, and G.~Salut,
``Maslov optimization theory: stochastic interpretation, particle
resolution,''
in: {\it Proceedings of the 11th Conf.
on Analysis and Optimization of
Systems: Discrete Event Systems},
vol.~199,
Springer Lecture Notes
in Control and Information Sciences,
1994,
pp.~312--318.

\item{32.}
M.~Akian, J.-P.~Quadrat, and M.~Voit,
``Bellman processes,''
in: {\it Proceedings of the 11th Conf.
on Analysis and Optimization of
Systems: Discrete Event Systems},
vol.~199,
Springer Lecture Notes
in Control and Information Sciences,
1994.

\item{33.}
J.~P.~Quadrat and Max-Plus working group,
``Max-plus algebra and applications
to system theory and optimal
control,''
in: {\it Proceedings of the International Congress of Mathematicians},
Z\"urich,
1994.

\item{34.}
M.~Akian, J.-P.~Quadrat, and M.~Voit,
``Duality between probability and optimization,''
in: {\it Idempotency},
Publ. of the Newton Institute,
Cambridge Univ. Press,
Cambridge,
1998,
pp.~331--353.

\item{35.}
P.~Del~Moral,
``Maslov optimization theory. Topological aspects,''
in: {\it Idempotency},
Publ. of the Newton Institute,
Cambridge Univ. Press,
Cambridge,
1998,
pp.~354--382.

\item{36.}
M.~Hasse,
``\"Uber die Behandlung graphentheoretischer Probleme unter
Verwendung der
Matrizenrechnung,''
{\it Wiss. Z. (Techn. Univer. Dresden)},
{\bf 10}
(1961),
1313--1316.

\item{37.}
J.~W.~Leech,
{\it Classical Mechanics},
Methuen's Monographs on Physical Subjects,
Methuen,
London, and
John Wiley and Sons,
New York,
1961.

\item{38.}
R.~P.~Feynman and A.~R.Hibbs,
{\it Quantum Mechanics and Path Integrals},
McGraw-Hill Book Company,
New~York,
1965.

\item{39.}
G.~W.~Mackey,
{\it The Mathematical Foundations of Quantum Mechanics. A Lecture-Note
Volume},
W.\ A.\ Benjamin,
New~York and Amsterdam,
1963.

\item{40.}
E.~Nelson,
{\it Probability Theory and Euclidian Field Theory},
Constructive Quantum Field Theory,
vol.~25,
Springer-Verlag,
Berlin,
1973.

\item{41.}
A.~N.~Vasil{\mz}ev,
{\it Functional Methods in Quantum Field Theory and Statistics\/}
[in Russian],
LGU,
Lenin\-grad,
1976.

\item{42.}
S.~Samborskii,
``The Lagrange problem
from the point of view of idempotent analysis,''
in: {\it Idempotency},
Publ. of the Newton Institute,
Cambridge Univ. Press,
Cambridge,
1998,
pp.~303--321.

\item{43.}
V.~P.~Maslov and S.~N.~Samborskii,
``Boundary-value problems for the stationary Hamilton--Jacobi and
Bellman equations,''
{\it Ukrain. Mat. Zh.} [{\it Ukrainian Math.~J.}],
{\bf 49}
(1997),
no.~3,
433--447.

\item{44.}
V.~P.~Maslov,
``A new approach to generalized solutions of nonlinear systems,''
{\it Dokl. Akad. Nauk SSSR\/} [{\it Soviet Math. Dokl.}],
{\bf 292}
(1987),
no.~1,
37--41.

\item{45.}
B.~Banaschewski and E.~Nelson,
``Tensor products and bimorphisms,''
{\it Canad. Math. Bull.},
{\bf 19}
(1976),
no.~4,
385--402.

\item{46.}
N.~Bourbaki,
{\it Topologie G{\'e}n{\'e}rale.
El{\'e}ments de Math{\'e}matique},
1{\`e}re partie, Livre~III,
Hermann,
Paris,
1960.

\item{47.}
V.~M.~Tikhomirov,
``Convex analysis,''
in: {\it Current Problems in Mathematics. Fundamental Directions\/}
[in Russian],
vol.~14,
Itogi Nauki i Tekhniki, VINITI,
Moscow,
1987,
pp.~5--101.

\item{48.}
J.~R.~Giles,
{\it Convex Analysis
with Application
in the Differentiation of Convex
Functions},
Pitman Publ.,
Boston--London--Melbourne,
1982.

\item{49.}
J.~Gunawardena,
``An introduction
to idempotency,''
in: {\it Idempotency},
Publ. of the Newton Institute,
Cambridge Univ. Press,
Cambridge,
1998,
pp.~1--49.

\item{50.}
G.~A.~Sarymsakov,
{\it Topological semifields\/}
[in Russian],
Fan,
Tashkent,
1969.

\item{51.}
S.~Eilenberg,
{\it Automata, Languages and Machines},
Academic Press,
1974.

\item{52.}
J.~H.~Conway,
{\it Regular Algebra and Finite Machines},
Chapman and Hall Math. Ser.,
London,
1971.

\item{53.}
K.~I.~Rosenthal,
{\it Quantales and Their Applications},
Pitman Research Notes
in Math. Ser,
Longman Sci.~\& Tech.,
1990.

\item{54.}
A.~Joyal and M.~Tierney,
``An extension of the Galois theory of Grothendieck,''
{\it Mem. Amer. Math. Soc.},
{\bf 51}
(1984),
no.~309,
VII.

\item{55.}
S.~N.~Pandit,
``A new matrix calculus,''
{\it SIAM J. Appl. Math.},
{\bf 9}
(1961),
no.~4,
632--639.

\item{56.}
N.~N.~Vorob{\mz}ev,
``The extremal matrix algebra,''
{\it Dokl. Akad. Nauk SSSR\/} [{\it Soviet Math. Dokl.}],
{\bf 152}
(1963),
no.~1,
24--27.

\item{57.}
N.~N.~Vorob{\mz}ev,
``The extremal algebra of positive matrices''
[in Russian],
{\it Elektronische Informationsverarbeitung und Kybernetik},
{\bf 3}
(1967),
39--71.

\item{58.}
N.~N.~Vorob{\mz}ev,
``The extremal algebra of nonnegative matrices''
[in Russian],
{\it Elektronische Informationsverarbeitung und Kybernetik},
{\bf 6}
(1970),
302--312.

\item{59.}
A.~A.~Korbut,
``Extremal spaces,''
{\it Dokl. Akad. Nauk SSSR\/} [{\it Soviet Math. Dokl.}],
{\bf 164}
(1965),
no.~6,
1229--1231.

\item{60.}
K.~Zimmermann,
``A general separation theorem
in extremal algebras,''
{\it Ekonomicko-Matematicky Obzor (Prague)},
{\bf 13}
(1977),
no.~2,
179--201.

\item{61.}
G.~Cohen, S.~Gaubert, and J.-P.~Quadrat,
``Hahn--Banach separation theorem
for Max-Plus semimodules,''
to appear in: {\it Proceedings of ``Conf\'erence en l'honneur du
Professeur Alain Bensoussan \`a l'occasion de son~60\`eme
anniversaire},
4--5 dec. 2000, Paris,
to be published
by IOS, the Netherlands.

\item{62.}
P.~Butkovi{\v c},
``Strong regularity of matrices--A survey of results,''
{\it Discrete Applied Math.},
{\bf 48}
(1994),
45--68.

\item{63.}
V.~Kolokoltsov,
``Idempotent structures in optimization,''
in: {\it Surveys in Modern Mathematics and Applications},
vol.~65,
{\it Proceedings of L.~S.~Pontryagin Conference,
Optimal Control},
Itogi Nauki i Tehniki,
VINITI,
Moscow,
1999,
pp.~118--174.

\item{64.}
R.~Bellman and W.~Karush,
``On a new functional transform
in analysis: the maximum transform,''
{\it Bull. Amer. Math. Soc.},
{\bf 67}
(1961),
501--503.

\item{65.}
I.~V.~Romanovskii,
``Optimization of the stationary control for a discrete deterministic
process''
[in Russian],
{\it Kibernetika}
(1967),
no.~2,
66--78.

\item{66.}
I.~V.~Romanovskii,
``Asymptotic behavior of a discrete deterministic process with
a continuous state space''
[in Russian],
{\it Optimal planning,
Trudy IM~SO Akademii Nauk SSSR},
{\bf 8}
(1967),
171--193.

\item{67.}
E.~Schr\"odinger,
``Quantization as an eigenvalue problem,''
{\it Annalen der Physik},
{\bf 364}
(1926),
361--376.

\item{68.}
E.~Hopf,
``The partial differential equation $u_t+uu_x=\mu u_{xx}$,''
{\it Comm. Pure Appl. Math.},
{\bf 3}
(1950),
201--230.

\item{69.}
P.-L.~Lions,
{\it Generalized Solutions of Hamilton--Jacobi Equations},
Pitman,
Boston,
1982.

\item{70.}
M.~G.~Crandall, H.~Ishii, and P.-L.~Lions,
``A user's guide
to viscosity solutions,''
{\it Bull. Amer. Math. Soc.},
{\bf 27}
(1992),
1--67.

\item{71.}
W.~H.~Fleming and H.~M.~Soner,
{\it Controlled Markov Processes and Viscosity Solutions},
Springer,
New York,
1993.

\item{72.}
A.~I.~Subbotin,
``Minimax solutions of first-order partial differential equations,''
{\it Uspekhi Mat. Nauk\/} [{\it Russian Math. Surveys\/}],
{\bf 51}
(1996),
no.~2,
105--138.

\item{73.}
{\sl Viscosity Solutions and Applications}
(I.~Capuzzo Dolcetta, P.-L.~Lions, editors),
vol.~1660,
Lecture Notes
in Math,
Springer,
1997.

\item{74.}
P.~D.~Lax,
``Hyperbolic systems of conservation laws II,''
{\it Comm. Pure Appl. Math.},
{\bf 10}
(1957),
537--566.

\item{75.}
O.~A.~Oleinik,
``Discontinuous solutions of nonlinear differential equations''
[in Russian],
{\it Uspekhi Mat. Nauk\/},
{\bf 12}
(1957),
no.~3,
3--73.

\item{76.}
O.~Viro,
``Dequantization of real algebraic geometry
on a logarithmic paper,''
3rd European Congress of Mathematics,
E-print {\tt math.AG/\allowlinebreak0005163}.

\item{77.}
G.~L.~Litvinov, V.~P.~Maslov, and G.~B.~Shpiz,
``Linear functionals on idempotent spaces. Algebraic approach,''
{\it Dokl. Ross. Akad. Nauk\/} [{\it Russian Acad. Sci. Dokl. Math.}],
{\bf 363}
(1998),
no.~3,
298--300.

\item{78.}
G.~L.~Litvinov, V.~P.~Maslov, and G.~B.~Shpiz,
{\it Idempotent Functional Analysis: An Algebraic Approach\/}
[in Russian],
International Sophus Lie Center,
Moscow,
1998.

\item{79.}
G.~L.~Litvinov, V.~P.~Maslov, and G.~B.~Shpiz,
{\it Idempotent Functional Analysis: An Algebraic Approach},
E-print {\tt math.FA/\allowlinebreak0009128},
2000.

\item{80.}
G.~L.~Litvinov, V.~P.~Maslov, and G.~B.~Shpiz,
``Tensor products of idempotent semimodules. An algebraic approach,''
{\it Mat. Zametki\/} [{\it Math. Notes\/}],
{\bf 65}
(1999),
no.~4,
572--585.

\item{81.}
A.~Grothendieck,
{\it Produits tensoriels topologiques et espaces nucl\'eaires},
vol.~16,
Mem. Amer. Math. Soc,
Amer. Math. Soc.,
Providence,
1955.

\endRefs

\authoraddress
{(V.~P.~Maslov) M.~V.~Lomonosov Moscow State University}

\authoraddress
{(G.~L.~Litvinov, G.~B.~Shpiz) International Sophus Lie Center}


\enddocument